\documentclass{amsart}
\usepackage{graphicx}
\usepackage{pdfsync}
\usepackage{fncylab}
\usepackage{footmisc}
\usepackage{mathrsfs}
\usepackage{amsthm}
\usepackage{amsmath}
\usepackage{amssymb}
\usepackage[T1]{fontenc}
    \newtheoremstyle{upright}%
        {1pt plus1pt minus1pt}%
        {1pt plus1pt minus1pt}%
        {\upshape}%
        {}%
        {\bfseries\scshape}%
        {\textbf{.}}%
        {1em}%
        {}%

\newtheorem{theorem}{Theorem}[section]

\newtheorem{lemma}[theorem]{Lemma}
\newtheorem{question}[theorem]{Question}

\newtheorem{proposition}[theorem]{Proposition}

\theoremstyle{definition}
\newtheorem{definition}[theorem]{Definition}
\newtheorem{example}[theorem]{Example}
\newtheorem{notation}[theorem]{Notation}
\newtheorem{remark}[theorem]{Remark}

\newcommand{\N}{\mathbb{N}} 
\newcommand{\R}{\mathbb{R}} 
\newcommand{\Z}{\mathbb{Z}} 

\newcommand{\cantor}{\mathscr{C}} 

\newcommand{\ofraci}[1]{\mathscr{O}_{#1}(\xio{#1},\theta_{#1}^0)}
\newcommand{\omegafraci}[1]{\Omega(F_{#1})}
\newcommand{\ofracixang}[3]{\mathscr{O}_{#1}(#2,#3)}
\newcommand{\ofraciang}[2]{\mathscr{O}_{#1}(\xio{#1},#2)}
\newcommand{\fraci}[1]{F_{#1}}
\newcommand{\omegafrac}{\Omega(F)}
\newcommand{\fprintfraci}[1]{\mathcal{F}_{#1}(\xio{#1},\theta_{#1}^0)}
\newcommand{\fprintfraciang}[2]{\mathcal{F}_{#1}(\xio{#1},#2)}

\newcommand{\kc}{K\!C}
\newcommand{\kci}[1]{K\!C_{#1}}
\newcommand{\ks}{K\!S}

\newcommand{\sksi}[1]{\mathcal{S}(\ks_{#1})}

\newcommand{\omegaksi}[1]{\Omega(\ks_{#1})}
\newcommand{\omegaks}{\Omega(\ks)}
\newcommand{\ksi}[1]{\ks_{#1}}
\newcommand{\xoo}{x_0^0}

\newcommand{\xio}[1]{x_{#1}^0}

\newcommand{\xii}[2]{x_{#1}^{#2}}

\newcommand{\fprintksi}[1]{\mathcal{F}_{#1}(\xio{#1},\theta_{#1}^0)}

\newcommand{\compseq}{\{\mathscr{O}_n(\xio{n},\theta^0)\}_{n=0}^\infty}
\newcommand{\compseqang}[1]{\{\mathscr{O}_n(\xio{n},#1)\}_{n=0}^\infty}
\newcommand{\compseqi}[1]{\{\mathscr{O}_n(\xio{n},\theta^0)\}_{n=#1}^\infty}
\newcommand{\compseqiang}[2]{\{\mathscr{O}_n(\xio{n},#2)\}_{n=#1}^\infty}
\newcommand{\compseqixang}[3]{\{\mathscr{O}_n(#2,#3)\}_{n=#1}^\infty}

\newcommand{\compseqsi}[1]{\{\mathscr{O}_n(\xio{n},\alpha^0)\}_{n=#1}^\infty}

\newcommand{\tern}[2]{[#1,#2]}

\newcommand{\mbfa}{\mathbf{a}}

\newcommand{\Sa}[1]{S_{\mathbf{#1}}} 
\newcommand{\celli}[2]{C_{#1,#2}} 

\newcommand{\alphaio}[1]{\alpha_{#1}^0}
\newcommand{\ssai}[1]{\mathcal{S}(S_{a,#1})}

\newcommand{\sai}[1]{S_{a,#1}}
\newcommand{\omegasa}[1]{\Omega(S_{#1})}
\newcommand{\omegasai}[2]{\Omega(S_{#1,#2})}
\newcommand{\omegasi}[1]{\Omega(S_{a,#1})}
\newcommand{\omegas}{\Omega(S_a)}
\newcommand{\osi}[1]{\mathscr{O}_{#1}(\xio{#1},\alpha_{#1}^0)}

\newcommand{\slopesa}[1]{\text{Slope}(S_{#1})}

\newcommand{\tfractal}{\mathscr{T}}
\newcommand{\tfraci}[1]{\mathscr{T}_{#1}}
\newcommand{\omegat}{\Omega(\mathscr{T})}
\newcommand{\omegati}[1]{\Omega(\mathscr{T}_{#1})}
\newcommand{\stfraci}[1]{\mathcal{S}(\tfraci{#1})}

\begin{document}
\title[Fractal Billiards]{The current state of fractal billiards}

\author[M. L. Lapidus]{Michel L. Lapidus}
\address{University of California, Department of Mathematics, 900 Big Springs Rd., Riverside, CA  92521-0135, USA}
\email{lapidus@math.ucr.edu}
\thanks{The work of M. L. Lapidus was partially supported by the National Science Foundation under the research grants DMS-0707524 and DMS-1107750, as well as by the Institut des Hautes Etudes Scientifiques (IHES) in Bures-sur-Yvette, France, where he was a visiting professor while this paper was written. The work of R. G. Niemeyer was partially supported by the National Science Foundation under the MCTP grant DMS-1148801, while a postdoctoral fellow at the University of New Mexico, Albuquerque.}
\author[R. G. Niemeyer]{Robert G. Niemeyer}
\address{University of New Mexico, Department of Mathematics \& Statistics, 311 Terrace NE, Albuquerque, NM  87131-0001, USA}
\email{niemeyer@math.unm.edu}

\keywords{fractal billiard, polygonal billiard, rational (polygonal) billiard, law of reflection, unfolding process, flat surface, translation surface, geodesic flow, billiard flow, iterated function system and attractor, self-similar set, fractal, prefractal approximations, Koch snowflake billiard, $T$-fractal billiard, self-similar Sierpinski carpet billiard, prefractal rational billiard approximations, sequence of compatible orbits, hook orbits, (eventually) constant sequences of compatible orbits, footprints, Cantor points, smooth points, elusive points, periodic orbits, periodic vs. dense orbits.}
\subjclass[2010]{Primary: 28A80, 37D40, 37D50, Secondary: 28A75, 37C27, 37E35, 37F40, 58J99.}

\begin{abstract}
If $D$ is a rational polygon, then the associated rational billiard table is given by $\Omega(D)$.  Such a billiard table is well understood.  If $F$ is a closed fractal curve approximated by a sequence of rational polygons, then the corresponding fractal billiard table is denoted by $\omegafrac$.  In this paper, we survey many of the results from [\textbf{LapNie1-3}] for the Koch snowflake fractal billiard $\omegaks$ and announce new results on two other fractal billiard tables, namely, the $T$-fractal billiard table $\omegat$ (see \cite{LapNie6}) and a self-similar Sierpinski carpet billiard table $\omegasa{a}$ (see \cite{CheNie}).

We build a general framework within which to analyze what we call a sequence of compatible orbits.  Properties of particular sequences of compatible orbits are discussed for each prefractal billiard $\omegaksi{n}$, $\omegati{n}$ and $\omegasai{a}{n}$, for $n=0,1,2\cdots$.  In each case, we are able to determine a particular limiting behavior for an appropriately formulated sequence of compatible orbits.  Such a limit either constitutes what we call a nontrivial path of a fractal billiard table $\omegafrac$ or else a periodic orbit of $\omegafrac$ with finite period. In our examples, $F$ will be either $\ks$, $\mathscr{T}$ or $S_{a}$.  Several of the results and examples discussed in this paper are presented for the first time.

We then close with a brief discussion of open problems and directions for further research in the emerging field of fractal billiards.

\end{abstract}

\maketitle
\setcounter{tocdepth}{2}
\tableofcontents

\section{Introduction}
\label{sec:Introduction}
This paper constitutes a survey of a collection of results from \cite{LapNie1,LapNie2,LapNie3} as well as the announcement of new results on the $T$-fractal billiard table $\omegat$ (see \cite{LapNie6}) and a self-similar Sierpinski carpet billiard table $\omegasa{a}$ (see \cite{CheNie}).

In \S\S\ref{sec:RationalBilliards} and \ref{sec:fractalGeometry}, we survey the necessary background material for understanding the remainder of the article.  More specifically, in \S\ref{sec:RationalBilliards}, we introduce the notion of a rational polygonal billiard, a translation surface determined from a rational polygonal billiard and discuss the consequence of a dynamical equivalence between the billiard flow and the geodesic flow.\footnote{The references \cite{GaStVo,Gut1,MasTa,Sm,Ta1,Ta2,Vo,Zo} provide an excellent survey of the various topics in the field of mathematical billiards, as well as specific results pertinent to the theory of rational polygonal billiards and associated translation surfaces or flat surfaces.} This dynamical equivalence allows us to express an orbit of a rational billiard table as a geodesic on an associated translation surface, and vice-versa, with the added benefit of being able to determine the reflection in certain types of vertices of a rational billiard table.  Furthermore, in \S\ref{sec:fractalGeometry}, we provide additional background material from the subject of fractal geometry necessary for understanding the construction of the Koch snowflake $\ks$, $T$-fractal $\mathscr{T}$,\footnote{The $T$-fractal $\mathscr{T}$ was previously studied in a different context in \cite{AcST}.} and a Sierpinski carpet $S_{\mbfa}$, as well as particular orbits and \textit{nontrivial paths}.

We then combine the background material presented in \S\S\ref{sec:RationalBilliards} and \ref{sec:fractalGeometry} to analyze the prefractal billiard tables $\omegaksi{n}$, $\omegati{n}$ and $\omegasai{a}{n}$, for $n=0,1,2\cdots$.  We begin by providing a general language for prefractal billiards and subsequently focus on determining sufficient conditions for what we are calling a \textit{sequence of compatible periodic orbits}.  While \S\S\ref{subsec:ThePrefractalKochSnowflakeBilliard}--\ref{subsec:AprefractalSelfSimilarSierpinskiCarpetBilliard} contain specific results and specialized definitions, there is an over-arching theme that is more fully developed in \S\ref{sec:FractalBilliards}.

In addition to providing a general language within which to analyze a fractal billiard, we discuss in \S\S\ref{subsec:TheKochSnowflakeFractalBilliard}--\ref{subsec:aSelf-SimilarSierpinskiCarpetBilliard} how one can determine well-defined orbits of $\omegaks$, $\omegat$ and $\omegasa{a}$, as well as nontrivial paths of $\omegaks$ and $\omegat$ that connect two elusive points of each respective billiard. Relying on the main result of \cite{Du-CaTy}, the second author and Joe P. Chen have shown that it is possible to determine a periodic orbit of a self-similar Sierpinski carpet billiard $\omegasa{a}$; additional results and proofs are forthcoming in \cite{CheNie}, but a synopsis is provided in \S\S\ref{subsec:AprefractalSelfSimilarSierpinskiCarpetBilliard} and \ref{subsec:aSelf-SimilarSierpinskiCarpetBilliard}.

Many of the results in \S\S\ref{sec:PrefractalRationalBilliards} and \ref{sec:FractalBilliards} are being announced for the first time.  Specifically, \S\S\ref{subsec:TheTFractalPrefractalBilliard} and \ref{subsec:TheTFractalBilliard} contain new results on the prefractal $T$-fractal billiard $\omegati{n}$ and the $T$-fractal billiard $\omegat$ (see \cite{LapNie6}); \S\S\ref{subsec:AprefractalSelfSimilarSierpinskiCarpetBilliard} and \ref{subsec:aSelf-SimilarSierpinskiCarpetBilliard} contain new results for a prefractal Sierpinski carpet billiard $\omegasai{a}{n}$ and self-similar Sierpinski carpet billiard $\omegasa{a}$, where $a$ is the single underlying scaling ratio (see \cite{CheNie}).  As these sections constitute announcements of new results on the respective prefractal and fractal billiards, we will provide in future papers \cite{CheNie, LapNie4, LapNie5, LapNie6} detailed statements and proofs of the results given therein.  Given the nature of the subject of \textit{fractal billiards}, we will close with a discussion of open problems and possible directions for future work, some of which are to appear in \cite{CheNie} and \cite{LapNie4,LapNie5,LapNie6}.

\section{Rational billiards}
\label{sec:RationalBilliards}
In this section, we will survey the dynamical properties of a billiard ball as it traverses a region in the plane bounded by a closed and connected polygon.  In the latter part of this article, we will remove the stipulation that the boundary be a polygon and focus on billiard tables having boundaries that are fractal or containing subsets that are fractal (while still being simple, closed and connected curves in the plane).

Under ideal conditions, we know that a point mass making a perfectly elastic collision with a $C^1$ surface (or curve) will reflect at an angle which is equal to the angle of incidence, this being referred to as the \textit{law of reflection}.


Consider a compact region $\Omega(D)$ in the plane with simple, closed and connected boundary $D$.  Then, $\Omega(D)$ is called a \textit{planar billiard} when $D$ is smooth enough to allow the law of reflection to hold, off of a set of measure zero (where the measure is taken to be the arc length measure on $D$).  Though the law of reflection implicitly states that the angles of incidence and reflection be determined with respect to the normal to the line tangent at the basepoint, we adhere to the equivalent convention in the field of mathematical billiards that the vector describing the position and velocity of the billiard ball (which amounts to the position and angle, since we are assuming unit speed) be reflected in the tangent to the point of incidence.\footnote{This is equivalent to reflecting the incoming vector through the normal to the tangent at the point of collision in the boundary.  We continue with the convention established in the text, since it  is more convenient in the context of polygonal billiards.  Moreover, the fact that the  equivalence relation on the phase space is defined in terms of the convention we have adopted, necessitates us continuing with this convention; see \cite{Sm} for a formal discussion of the equivalence relation defined on the phase space $\Omega(D)\times S^1$.}  That is, employing such a law in order to determine the path on which the billiard ball departs after impact essentially amounts to identifying certain vectors. Such an equivalence relation is denoted by $\sim$ and, in the context of a polygonal billiard, is discussed below in more detail.

For the remainder of the article, unless otherwise indicated, when $D$ is a simple, closed, connected and piecewise smooth curve so as to allow the law of reflection to hold (off finitely many points), we assume $D$ is a closed and connected polygon. In such a case, we will refer to $\Omega(D)$ as a \textit{polygonal billiard}.

\begin{figure}
\begin{center}
\includegraphics[scale = .5]{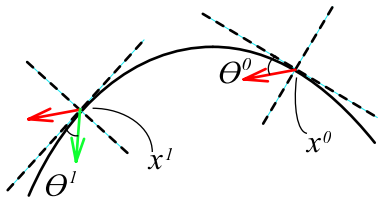}
\end{center}
\caption[Recovering the law of reflection]{A billiard ball traverses the interior of a billiard and collides with the boundary.  The velocity vector is pointed outward at the point of collision. The resulting direction of flow is found by either reflecting the vector through the tangent or by reflecting the incidence vector through the normal and reversing the direction of the vector.  We use the former method in this paper.}
\label{fig:billiardMap}
\end{figure}

One may express the law of reflection in terms of equivalence classes of vectors by identifying two particular vectors that form an equivalence class of vectors in the unit tangent bundle corresponding to the billiard table $\Omega(D)$; see Figure \ref{fig:billiardMap}. (See \cite{Sm} for a detailed discussion of this equivalence relation on the unit tangent bundle $\Omega(D)\times S^1$.)

Denote by $S^1$ the unit circle, which we let represent all the possible directions (or angles) in which a billiard ball may initially move.  To clearly understand how one forms equivalence classes from elements of $\Omega(D)\times S^1$, we let $(x,\theta),(y,\gamma)\in \Omega(D)\times S^1$ and say that $(x,\theta)\sim (y,\gamma)$ if and only if $x=y$ and one of the following is true:
\begin{enumerate}
\item{$x=y$ is not a vertex of the boundary $D$ and $\theta = \gamma$;}
\item{$x=y$ is not a vertex of the boundary $D$, but $x=y$ is a point on a segment $s_i$ of the polygon $D$ and $\theta = r_i(\gamma)$, where $r_i$ denotes reflection in the segment $s_i$;}
\item{If $x=y$ is a vertex of $D$, then we identify $(x,\theta)$ with $(y,g(\gamma))$ for every $g$ in the group generated by reflections in the two adjacent sides having $x$ (or $y$) as a common vertex.}
\end{enumerate}

\noindent For now, we shall denote by $[(x,\theta)]$ the equivalence class of $(x,\theta)$, relative to the equivalence relation $\sim$.

The collection of vertices of $\Omega(D)$ forms a set of zero measure (when we take our measure to be the arc-length measure on $D$), since there are finitely many vertices.



The phase space for the billiard dynamics is given by the quotient space $(\Omega(D)\times S^1)/\sim$.  In practice, one restricts his or her attention to the space $(D\times S^1)/\sim$.  The billiard flow on $(D\times S^1)/\sim$ is determined from the continuous flow on $\Omega(D)\times S^1)/\sim$ as follows.  Let $x^0$ be an initial basepoint, $\theta^0$ be an initial direction and $\varphi_t(x^0,\theta^0)$ be a flow line corresponding to these initial conditions in the phase space $(\Omega(D)\times S^1)/\sim$.    The values $t_j$ for which $\varphi_{t_j}(x^0,\theta^0)\in (D\times S^1)/\sim$ constitute the \textit{return times} (i.e., times at which $\varphi_t(x^0,\theta^0)$ returns to the section, or intersects it in a non-tangential way).  Then, the discrete map $f^{t_j}(x^0,\theta^0)$ constitutes the \textit{section map}.  In terms of the configuration space, $f^{t_j}(x^0,\theta^0)$ constitutes the point and angle of incidence in the boundary $D$.  Since $\Omega(D)$ is the billiard table and we are interested in determining the collision points, it is only fitting that such a map be called the \textit{billiard map}.  More succinctly, we denote $f^{t_j}$ by $f^j$ and, in general, such a map is called the \textit{Poincar\'e map} and the section is called the \textit{Poincar\'e section}.  Furthermore, the obvious benefit of having a visual representation of $f^j(x^0,\theta^0)$ in the configuration space is exactly why one restricts his or her attention to the section $(D\times S^1)/\sim$.  Specifically, all one really cares about in the end, from the perspective of studying a planar billiard, are the collision points, which are clearly determined by the billiard map.

In order to understand how one determines the next collision point and direction of travel, we must further discuss the \textit{billiard map} $f_D$. As previously discussed, $f_D:(D\times S^1)/\sim\,\to (D\times S^1)/\sim$, where the equivalence relation $\sim$ is the one introduced above.  More precisely, if $\theta^0$ is an inward pointing vector at a basepoint $x^0$, then $(x^0,\theta^0)$ is the representative element of the equivalence class $[(x^0,\theta^0)]$.  The billiard map then acts on $(D\times S^1)/\sim$ by mapping $[(x^k,\theta^k)]$ to $[(x^{k+1},\theta^{k+1})]$, where $x^k$ and $x^{k+1}$ are collinear in the direction determined by $\theta^k$ and where $\theta^{k+1}$ is the reflection of angle $\theta^k$ through the tangent at $x^{k+1}$.  In general, we have $f_D^k[(x^0,\theta^0)] = [(x^k,\theta^k)]$, for every $k\geq 0$.

\begin{remark}
In the sequel, we will simply refer to an element $[(x^k,\theta^{k})]\in (\Omega(D)\times S^1)/\sim$ by $(x^k,\theta^k)$, since the vector corresponding to $\theta^k$ is inward pointing at the basepoint $x^k$.  So as not to introduce unnecessary  notation, when we discuss the billiard map $f_{\fraci{n}}$ corresponding to the $n$th prefractal billiard $\omegafraci{n}$ approximating a fractal billiard $\omegafrac$, we will simply write $f_{\fraci{n}}$ as $f_n$.  When discussing the discrete billiard flow on $(\omegafraci{n}\times S^1)/\sim$, the $k$th point in an orbit $(x^k,\theta^k)\in (\omegafraci{n}\times S^1)/\sim$ will instead be denoted by $(\xii{n}{k_n},\theta_n^{k_n})$, in order to keep track of the space such a point belongs to (namely, with our present convention, $(\omegafraci{n}\times S^1)/\sim$).  Specifically, $k_n$ refers to the number of iterates of the billiard map $f_n$ necessary to produce the pair $(\xii{n}{k_n},\theta_n^{k_n})$. An initial condition of an orbit of $\omegafraci{n}$ will always be referred to as $(\xio{n},\theta_n^0)$.
\end{remark}

In what follows, we are presupposing an orbit can be formed by iterating the billiard map forward in time and backwards in time, whenever $f_n^{-k}(\xio{n},\theta_n^0)$ is defined.

An orbit making finitely many collisions in the boundary is called a \textit{closed orbit}.  If, in addition, there exists $m\in\mathbb{Z}$ such that $f_D^m(x^0,\theta^0) = (x^0,\theta^0)$, then the resulting orbit is called \textit{periodic}; the smallest positive integer $m$ such that $f^m_D(x^0,\theta^0)=(x^0,\theta^0)$ is called the \textit{period} of the periodic orbit. In the event that a basepoint $x^j$ of $f_D^j(x^0,\theta^0)$ is a corner of $\Omega(D)$ (that is, a vertex of the polygonal boundary $D$) and reflection cannot be determined in a well-defined manner, then the resulting orbit is said to be \textit{singular}.  In addition, if there exists a positive integer $k$ such that the basepoint $x^{-k}$ of $f_D^{-k}(x^0,\theta^0)$ is a corner of $\Omega(D)$ (here, $f_D^{-k}$ denotes the $k$th inverse iterate of $f_D$), then the resulting orbit is closed and the path traced out by the billiard ball connecting $x^j$ and $x^{-k}$ is called a \textit{saddle connection}.  Finally, we note that a periodic orbit with period $m$ is a closed orbit for which reflection is well defined at each basepoint $x^i$ of $f_D^i(x^0,\theta^0)$, $0\leq i\leq m$ and $f_D^m(x^0,\theta^0) = (x^0,\theta^0)$.

We say that an orbit $\mathscr{O}(x^0,\theta^0)$ is \textit{dense} in a rational billiard table $\Omega(D)$ if the path traversed (forward and backward in time) by the billiard ball in $\Omega(D)$ is dense in $\Omega(D)$. That is, the closure of the set of points comprising the path traversed by the billiard ball is exactly $\Omega(D)$.  Likewise, the points of incidence (i.e., the footprint) of a dense orbit will be dense in the boundary $D$, as explained in Remark \ref{rmk:DenseOrbits}.

\begin{remark}
\label{rmk:DenseOrbits}
Consider a rational polygonal billiard $\Omega(D)$.  The associated translation surface $\mathcal{S}(D)$ can be constructed as described in \S\ref{subsec:translationStructuresandTranslationSurfaces}.  As we will show in \S\ref{subsec:UnfoldingABilliardOrbit}, the geodesic flow on a translation surface is dynamically equivalent to the billiard flow.  A dense orbit will have an initial direction preventing the path from being parallel to any side of $\Omega(D)$ (except, possibly, for finitely many initial directions, and hence, for a measure-zero set).  The corresponding path on the associated translation surface\footnote{See \S\ref{subsec:translationStructuresandTranslationSurfaces} for an explanation of what constitutes a translation surface.} must also be dense in the surface.  Since the path on the surface is arbitrarily close to every side appropriately identified with another side of a copy of $\Omega(D)$ and not parallel to any side, the path will be transversal with respect to each side.  Thus, the collection of basepoints of a dense orbit must be dense in $D$.
\end{remark}

\begin{definition}[Footprint of an orbit]
\label{def:footprintOfAnOrbit}
Let $\mathscr{O}_D(x^0,\theta^0)$ be an orbit of a billiard $\Omega(D)$ with an initial condition $(x^0,\theta^0)\in D\times S^1$.  Then the trace of an orbit on the boundary $D$,
\begin{equation}
\mathscr{O}_D(x^0,\theta^0)\cap D,
\end{equation}
is called the \textit{footprint} of the orbit $\mathscr{O}_D(x^0,\theta^0)$ and is denoted by $\mathcal{F}_D(x^0,\theta^0)$.  When we are only interested in a prefractal billiard $\Omega(F_n)$, we denote the footprint of an orbit by $\mathcal{F}_n(\xio{n},\theta^0_n)$.
\end{definition}

For the remainder of the article, when discussing polygonal billiards, we will focus our attention on what are called \textit{rational polygonal billiards}, or, more succinctly, \textit{rational billiards}.

\begin{definition}[Rational polygon and rational billiard]
If $D$ is a nontrivial connected polygon such that for each interior angle $\theta_j$ of $D$ there are relatively prime integers $p_j \geq 1$ and $q_j\geq 1$ such that $\theta_j = \frac{p_j}{q_j} \pi$, then we call $D$ a \textit{rational polygon} and $\Omega(D)$ a \textit{rational billiard}.
\label{def:ratBilliard}
\end{definition}

\subsection{Translation surfaces and properties of the flow}
\label{subsec:translationStructuresandTranslationSurfaces}
In this subsection, we will discuss what constitutes a translation surface and how to construct a translation surface from a rational billiard.  Then, in \S\ref{subsec:UnfoldingABilliardOrbit}, we will see how to relate the continuous billiard flow on $(\Omega(D)\times S^1)/\sim$ with the geodesic flow on the associated translation surface.

\begin{definition}[Translation structure and translation surface]
\label{def:translationStructure}
Let $M$ be a compact, connected, orientable surface.  A \textit{translation structure} on $M$ is an atlas $\omega$, consisting of charts of the form $(U_\alpha,\varphi_\alpha)_{\alpha\in\mathscr{A}}$, where $U_\alpha$ is a domain (i.e., a connected open set) in $M$ and $\varphi_\alpha$ is a homeomorphism from $U_\alpha$ to a domain in $\mathbb{R}^2$, such that the following conditions hold:

\begin{enumerate}
\item{The collection $\{U_\alpha\}_{\alpha\in\mathscr{A}}$ covers the whole surface $M$ except for finitely many points $z_1,z_2,...,z_k$, called \textit{singular points};}
\item{all coordinate changing functions are translations in $\mathbb{R}^2$;}
\item{the atlas $\omega$ is maximal with respect to properties $(1)$ and $(2)$;}
\item{for each singular point $z_j$, there is a positive integer $m_j$, a punctured neighborhood $\dot{U}_j$ of $z_j$ not containing other singular points, and a map $\psi_j$ from this neighborhood to a punctured neighborhood $\dot{V}_j$ of a point in $\mathbb{R}^2$ that is a shift in the local coordinates from $\omega$, and is such that each point in $\dot{V}_j$ has exactly $m_j$ preimages under $\psi_j$.}
\end{enumerate}

\noindent We say that a connected, compact surface equipped with a translation structure is a \textit{translation surface}.
\end{definition}

\begin{remark}
Note that in the literature on billiards and dynamical systems, the terminology and definitions pertaining to this topic are not completely uniform; see, for example, \cite{GaStVo,Gut1,GutJu1,GutJu2,HuSc,Mas,MasTa,Ve1,Ve2,Vo,Zo}. (We note that in \cite{MasTa} and \cite{Zo}, `translation surfaces' are referred to as `flat surfaces'.) We have adopted the above definition for clarity and the reader's convenience.
\end{remark}

We now discuss how to construct a translation surface from a rational billiard.  Consider a rational polygonal billiard $\Omega(D)$ with $k$ sides and interior angles $\frac{p_j}{q_j}\pi$ at each vertex $z_j$, for $1\leq j\leq k$, where the positive integers $p_j$ and $q_j$ are relatively prime.  The linear parts of the planar symmetries generated by reflection in the sides of the polygonal billiard $\Omega(D)$ generate a  dihedral group $\mathscr{D}_N$, where $N:=\text{lcm} \{q_j\}_{j=1}^k$ (the least common multiple of the $q_j$'s). Next, we consider $\Omega(D)\times \mathscr{D}_N$ (equipped with the product topology).  We want to glue `sides' of $\Omega(D)\times  \mathscr{D}_N$ together and construct a natural atlas on the resulting surface $M$ so that $M$ becomes a translation surface.

\begin{figure}
\begin{center}
\includegraphics[scale=.5]{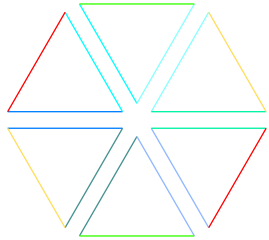}
\end{center}
\caption{The equilateral triangle billiard $\Omega(\Delta)$ can be acted on by a particular group of symmetries to produce a translation surface that is topologically equivalent to the flat torus. In this figure, we see that opposite and parallel sides are identified in such a way that the orientation is preserved.  This allows us to examine the geodesic flow on the surface.  We will see in \S\ref{subsec:UnfoldingABilliardOrbit} that the geodesic flow on the translation surface is dynamically equivalent to the continuous billiard flow.}
\label{fig:sixEquiSurface}
\end{figure}

As a result of the identification, the points of $M$ that correspond to the vertices of $\Omega(D)$ constitute (removable or nonremovable) conic singularities of the surface.  Heuristically, $\Omega(D)\times \mathscr{D}_N$ can be represented as $\{r_j\Omega(D)\}_{j=1}^{2N}$, in which case it is easy to see what sides are made equivalent under the action of $\sim$.  That is, $\sim$ identifies opposite and parallel sides in a manner which preserves the orientation.  See Example \ref{exa:EquilateralTriangleFlatSurface} and Figure \ref{fig:sixEquiSurface} for an example of a translation surface constructed from the equilateral triangle billiard $\Omega(\Delta)$.

\begin{example}
\label{exa:EquilateralTriangleFlatSurface}
Consider the equilateral triangle $\Delta$.  The corresponding billiard is denoted by $\Omega(\Delta)$.  The interior angles are $\{\frac{\pi}{3},\frac{\pi}{3},\frac{\pi}{3}\}$.  Hence, the group acting on $\Omega(\Delta)$ to produce the translation surface is the dihedral group $\mathscr{D}_3$. The resulting translation surface is topologically equivalent to the flat torus. We will make use of this fact in the sequel.
\end{example}

\subsection{Unfolding a billiard orbit and equivalence of flows}
\label{subsec:UnfoldingABilliardOrbit}
\begin{figure}
\begin{center}
\includegraphics[scale=.35]{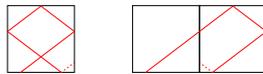}
\end{center}
\caption{Partially unfolding an orbit of the square billiard $\Omega(Q)$.  The `R' is shown so as to provide the reader with a frame of reference.}
\label{fig:PartiallyUnfoldingAnOrbitOfTheSquareBilliard}
\end{figure}

\begin{figure}
\begin{center}
\includegraphics[scale=.35]{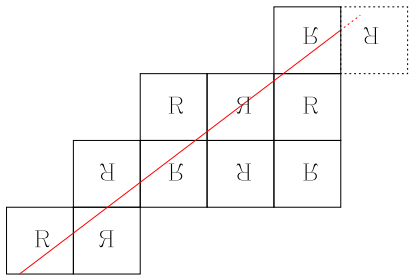}
\end{center}
\caption{Unfolding an orbit of the square billiard $\Omega(Q)$.}
\label{fig:UnfoldingAnOrbitOfTheSquareBilliard}
\end{figure}

Consider a rational polygonal billiard $\Omega(D)$ and an orbit $\mathscr{O}(x^0,\theta^0)$.  Reflecting the billiard $\Omega(D)$ and the orbit in the side of the billiard table containing the basepoint $x^1$ of the orbit (or an element of the footprint of the orbit) partially unfolds the orbit $\mathscr{O}(x^0,\theta^0)$; see Figure \ref{fig:PartiallyUnfoldingAnOrbitOfTheSquareBilliard} for the case of the square billiard $\Omega(Q)$. Continuing this process until the orbit is a straight line produces as many copies of the billiard table as there are elements of the footprint; see Figure \ref{fig:UnfoldingAnOrbitOfTheSquareBilliard}. That is, if the period of an orbit $\mathscr{O}(x^0,\theta^0)$ is some positive integer $p$, then the number of copies of the billiard table in the unfolding is also $p$. We refer to such a straight line as the \textit{unfolding of the billiard orbit}.

\begin{figure}
\begin{center}
\includegraphics[scale=.35]{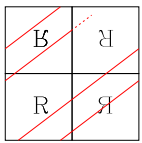}
\end{center}
\caption{Rearranging the unfolded copies of the unit square from Figure \ref{fig:UnfoldingAnOrbitOfTheSquareBilliard} and correctly identifying sides so as to recover the flat torus, we see that the unfolded orbit corresponds to a closed geodesic of the translation surface.}
\label{fig:equivalenceOfGeodesicFlow}
\end{figure}

Given that a rational billiard $\Omega(D)$ can be acted on by a dihedral group $\mathscr{D}_N$ to produce a translation surface in a way that is similar to unfolding the billiard table, we can quickly see how the billiard flow is dynamically equivalent to the geodesic flow; see Figure \ref{fig:equivalenceOfGeodesicFlow} and the corresponding caption.

One may modify the notion of ``reflecting'' so as to determine orbits of billiard tables tiled by a rational polygon $D$.  As an example, we consider the unit-square billiard table. An appropriately scaled copy of the unit-square billiard table can be tiled by the unit-square billiard table by making successive reflections in the sides of the unit square.  One may then unfold an orbit of the unit-square billiard table into a larger square billiard table. When the unfolded orbit of the original unit-square billiard intersects the boundary of the appropriately scaled (and larger) square, then one continues unfolding the billiard orbit in the direction determined by the law of reflection (that is, assuming the unfolded orbit is long enough to reach a side of the larger square).  We will refer to such an unfolding as a \textit{reflected-unfolding}.

\begin{figure}
\begin{center}
\includegraphics[scale=.35]{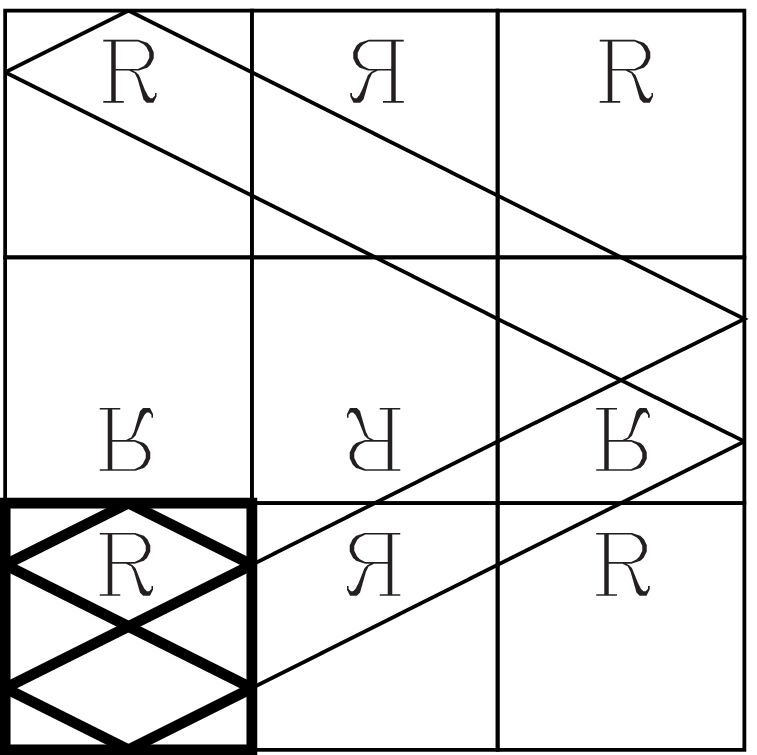}
\end{center}
\caption{Unfolding the orbit of the unit-square billiard in a (larger) scaled copy of the unit-square billiard. This constitutes an example of a reflected-unfolding. The edges of the original unit-square billiard table and the segments comprising the orbit have been thickened to provide the reader with a frame of reference.}
\label{fig:aReflectedUnfoldingInTheSquare}
\end{figure}

We may continue this process in order to form an orbit of a larger scaled square billiard table. Suppose that an orbit $\mathscr{O}(x^0,\theta^0)$ has period $p$.  The footprint of the orbit is then $\mathscr{F}_B(x^0,\theta^0)=\{f^i_B(x^0,\theta^0)\}_{i=0}^{p-1}$.   If $s$ is a positive integer (i.e., $s\in \mathbb{N}$), then the footprint $\mathscr{F}_B^s(x^0,\theta^0):=\{f^i_B(x^0,\theta^0)\}_{i=0}^{s(p-1)}$ of an orbit constitutes the footprint of an orbit that traverses the same path $s$-many times.  For sufficiently large $s\in \mathbb{N}$, an orbit that traverses the same path as an orbit $\mathscr{O}(x^0,\theta^0)$ $s$-many times can be reflected-unfolded in an appropriately scaled square billiard table to form an orbit of the larger billiard table; see Figure \ref{fig:aReflectedUnfoldingInTheSquare}.  Such a tool is useful in understanding the relationship between the billiard flow on a rational polygonal billiard $\Omega(D)$ and a billiard table tiled by $D$, and will be particularly useful in understanding the nature of particular orbits of a self-similar Sierpinski carpet billiard in \S\ref{subsec:aSelf-SimilarSierpinskiCarpetBilliard}.


\begin{figure}
\begin{center}
\includegraphics[scale=.6]{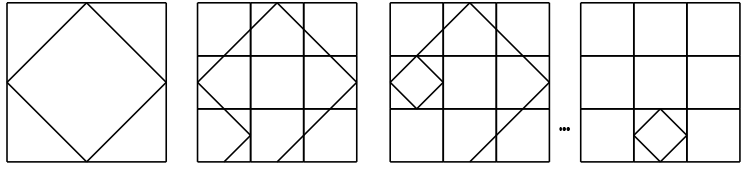}
\end{center}
\caption{Illustrated in this figure is the process of folding up an orbit of a square billiard table, as discussed at the end of \S\ref{subsec:UnfoldingABilliardOrbit}.  In the first image, we see an orbit of unit-square billiard table.  Partitioning the unit square into nine equally sized squares, we see that we can fold up the orbit by making successive reflections in the sides of the squares comprising the partition.  Using sufficiently many reflections results in an orbit of one of the squares of the partition.}
\label{fig:foldedUpOrbit}
\end{figure}

As one may expect, if $D$ is a rational polygon that tiles a billiard table $\Omega(R)$, then an orbit of $\Omega(R)$ may be folded up to form an orbit of $\Omega(D)$. This is done by making successive reflections in $D$, the result being an orbit of $\Omega(D)$; see Figure \ref{fig:foldedUpOrbit} for the case of a square billiard table.

\section{The fractals of interest}
\label{sec:fractalGeometry}
We are primarily interested in fractals with boundaries either partially or completely comprised of self-similar sets and fractals that are self-similar.  So as to make the material discussed in \S{\ref{sec:PrefractalRationalBilliards}--\ref{sec:FractalBilliards}} more accessible, we provide a few basic definitions from the subject of fractal geometry.

\begin{definition}
\label{def:Contraction}
Let $(X,d)$ be a metric space and $\phi:X\to X$.

\begin{itemize}
\item[(i)]{(Contraction)}.\quad
If there exists $0 < c < 1$ such that
\begin{align}
\notag d(\phi(x),\phi(y)) &\leq c d(x,y)
\end{align}
\noindent for every $x,y\in X$, then $\phi$ is called a  \textit{contraction} (or \textit{contraction mapping}).

\item[(ii)]{(Similarity contraction)}.\quad If there exists $0<c<1$ such that
\begin{align}
\notag d(\phi(x),\phi(y)) &= c d(x,y),
\end{align}
\noindent for every $x,y\in X$, then $\phi$ is called a \textit{similarity contraction}. This unique value $c\in (0,1)$ is called the \textit{scaling ratio} of $\phi$.
\end{itemize}
\end{definition}

\begin{definition}
\label{def:iteratedFunctionSystemAndSelfSimilarSet}
Let $(X,d)$ be a complete metric space.

\begin{itemize}
\item[(i)]{(Iterated function system and attractor)}.\quad Let $\{\phi_i\}_{i=1}^k$ be a family of contractions defined on $X$. Then $\{\phi_i\}_{i=1}^k$ is called an \textit{iterated function system} (IFS).

An iterated function system is so named because the map $\Phi:\mathbf{K}\to\mathbf{K}$, given by  $\Phi(\cdot):=\bigcup_{i=1}^k \phi_i(\cdot)$ and defined on the space $\mathbf{K}$ of nonempty compact subsets of $X$, can be composed with itself.  Indeed, for each $m\in \N$, we have
\begin{align}
\Phi^m(\cdot) &= \bigcup_{i_1=1}^k...\bigcup_{i_m=1}^k \phi_{i_1}\circ\cdots\circ\phi_{i_m}(\cdot).
\label{eqn:PhiMiterate}
\end{align}
Furthermore, there exists a unique nonempty compact set $F\subset X$ (i.e., $F\in \mathbf{K}$), called the \textit{attractor} of the IFS, such that
\begin{align}
F=\Phi(F)&:= \bigcup_{i=1}^k \phi_i(F).
\end{align}
\item[(ii)]{(Self-similar system and self-similar set)}.\quad In the special case where each $\phi_i$ is a contraction similarity, for $i=1,...,k$, then the IFS $\{\phi_i\}_{i=1}^k$ is said to be a \textit{self-similar system} and its attractor $F$ is called a \textit{self-similar set} (or a \textit{self-similar subset} of $X$).
\end{itemize}
\end{definition}

If $X$ is complete, then so is $\mathbf{K}$ (equipped with the Hausdorff metric\footnote{See \cite{Ba} for details on the Hausdorff metric.}) and hence, since it can be shown that $\Phi:\mathbf{K}\to \mathbf{K}$ is a contraction, it follows from the contraction mapping theorem that $\Phi$ has a unique fixed point (thereby justifying the definition of the attractor $F$ above) and that for any $E\in \mathbf{K}$, $\Phi^m(E)\to F$, as $m\to\infty$ (where, as in Equation (\ref{eqn:PhiMiterate}), $\Phi^m$ is the $m$th iterate of $\Phi$). (See \cite{Hut}.)

We state the next property in the special case which will be of interest to us, namely, that of an IFS in a Euclidean space.

\begin{theorem}[{\cite{Hut}; see also \cite[Thm. 9.1]{Fa}}]
\label{thm:FalconersTheorem}
Consider an iterated function system given by contractions $\{\phi_i\}_{i=1}^k$, each defined on a compact set $D\subseteq \mathbb{R}^n$, such that $\phi_i(D)\subseteq D$ for each $i\leq k$, and with attractor $F$. Then $F\subseteq D$ and in fact,
\begin{align}
F&= \bigcap_{m=0}^\infty \Phi^m(E)
\end{align}
\noindent for every set $E\in \mathbf{K}$ such that $\phi_i(E)\subseteq E$ for all $i\leq k$.  Here, the transformation $\Phi:\mathbf{K}\to\mathbf{K}$ is given as in part \emph{(}i\emph{)} of Definition \ref{def:iteratedFunctionSystemAndSelfSimilarSet}.
\end{theorem}

\begin{notation}
Suppose $F$ is a fractal set.  Then, the $n$th prefractal approximation of $F$ is denoted by $F_n$.  In the case of a self-similar fractal $F$, the $n$th \textit{prefractal approximation} of $F$ is usually defined by $\bigcap_{m=0}^n \Phi^m(E)$, where $E\in\mathbf{K}$.
\end{notation}

Not every fractal is self-similar or embedded in Euclidean space.  However, such sets represent an important collection of examples of fractal sets.  In the next subsection, we will discuss the fractal subsets (self-similar or not) of $\mathbb{R}$ or of $\mathbb{R}^2$ of direct interest to us in this paper.

\subsection{Cantor sets}
\label{subsec:CantorSets}
A Cantor set is a set with very rich and counter-intuitive properties; topologically, it is a compact and totally disconnected (i.e., perfect) space. In order to illustrate some of the properties that make a Cantor set so interesting, we refer to the canonical example of a Cantor set: the ternary Cantor set.  We focus on three methods for constructing the ternary Cantor set: 1) by tremas, 2) as the unique fixed point attractor of an iterated function system, and 3) in terms of an alphabet.

Before we discuss the ternary Cantor set, we mention that this set was first discovered by Henry J. S. Smith in 1875.  Later, in 1881, Vito Volterra independently rediscovered the ternary Cantor set.  Smith's and Volterra's records being obscured over the years for one reason or another, it was the German mathematician Georg Cantor whom, in 1883, history credits with the discovery of a bounded, totally disconnected, perfect and uncountable set with measure zero, that is now commonly referred to as ``the Cantor set''.

\begin{figure}
\begin{center}
\includegraphics[scale=.5]{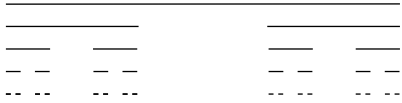}
\end{center}
\caption{The ternary Cantor set.}
\label{fig:TernaryCantorSet}
\end{figure}

We now proceed to construct the ternary Cantor set, hereafter denoted by $\mathscr{C}$, by the method known as \textit{construction by tremas}, which is Latin for `cuts'.  Begin with the unit interval $I$ and remove the middle open third $(\frac{1}{3},\frac{2}{3})$ from $I$, leaving the two closed intervals $[0,\frac{1}{3}]$ and $[\frac{2}{3},1]$.  Next, remove the middle open ninth from each closed subinterval.  What remains are the closed intervals $[0,\frac{1}{9}]$, $[\frac{2}{9},\frac{1}{3}]$, $[\frac{2}{3},\frac{7}{9}]$, $[\frac{8}{9},1]$.  Continuing this process ad infinitum, we construct the ternary Cantor set; see Figure \ref{fig:TernaryCantorSet}.

One may also construct $\cantor$ by utilizing an appropriately defined iterated function system.  Consider the following contraction maps defined on the real line $\R$:
\begin{align}
\phi_1(x)= \frac{1}{3} x, \quad & \quad \phi_2(x) = \frac{1}{3} x + \frac{2}{3}.
\end{align}
\noindent Then, $\lim_{n\to\infty} \Phi^n(I)= \mathscr{C}$, where $\Phi$ is given as in part (i) of Definition \ref{def:iteratedFunctionSystemAndSelfSimilarSet}.  Moreover, since $\{\phi_i\}_{i=1}^2$ is a family of similarity contractions and $\mathscr{C} = \Phi(\mathscr{C})$, we have that $\mathscr{C}$ is a self-similar set.

A third---and equivalent---construction of the ternary Cantor set can be given in terms of the symbols $l$, $c$, and $r$.  Recall that the elements of $\mathbb{R}$ can be expressed in terms of a base-$3$ number system.  We focus our attention on elements of the unit interval $I$.  So-called ternary numbers\footnote{An element $x\in I$ is a \textit{ternary number} if $x=\frac{p}{3^y}$, $0\leq p\leq 3^y$, $p,y\in \mathbb{N}$.} in $I$ have two equivalent expansions: one that is finite and one that is infinite.  For example, $\frac{1}{3}$ can be written in base-$3$ as $0.1$ or, equivalently, as $0.0\overline{2}$ (where the overbar indicates that the digit $2$ is repeated infinitely often).

We next discuss a similar \textit{addressing system} that has the benefit of preventing ternary numbers from having a finite representation.  The characters $l$, $c$ and $r$ are to remind the reader of choosing \textit{left}, \textit{center} and \textit{right}. We identify an element of the unit interval $I$ by an infinite address that indicates \textit{where} in $I$ the element is located.  Motivated by the construction of $\mathscr{C}$ by tremas, one can identify any point of $I$ by an infinite address consisting of the characters $l$, $c$ and $r$.  While elements of $\mathscr{C}$ can be represented by infinite addresses consisting of $c$'s, we make the stipulation that no element of $\mathscr{C}$ will be represented by an infinite address containing $c$'s.\footnote{In other words, we do not allow an element of $\cantor$ to be approximated by a sequence $\{z_i\}_{i=1}^\infty$ of elements of $\cantor^c$, where $\cantor^c$ is the complement of the ternary Cantor set in $I$.}  Moreover, this method of representing elements of $I$ (or $\cantor$) provides every element with an infinite representation and never a finite representation.

\begin{example}
The values $\frac{1}{4}$, $\frac{1}{3}$ and $\frac{1}{2}$ have the ternary representations $\overline{lr}$, $l\overline{r}$ and $\overline{c}$, respectively.\footnote{Equivalently, $\frac{1}{3}$ has a representation given by $c\overline{l}$. Although, we will not consider this as a representation for $\frac{1}{3}$ on account of $\frac{1}{3}\in \mathscr{C}$.}  While $\frac{1}{3}$ has a finite ternary expansion given by $0.2=0.0\overline{1}$, it does not have a finite ternary representation.  It should be noted that elements like $\frac{1}{4}$ and $\frac{1}{2}$ will play an important role in our analysis of the Koch snowflake fractal billiard.  The occurrence of infinitely many $c$'s or infinitely many $l$'s and $r$'s is critical to developing some of the theory regarding the Koch snowflake fractal billiard.
\end{example}

So that some of the results concerning the Koch snowflake fractal billiard can be more succinctly expressed, we introduce a notation used for describing a value's \textit{type of ternary representation}.

\begin{notation}[The type of ternary representation]
\label{not:typeOfTernaryRepresentation}
The \textit{type of ternary representation} can be defined as follows.  If $x\in I$, then the first coordinate of $[\cdot, \cdot]$ describes the characters that occur infinitely often and the second coordinate of $[\cdot,\cdot]$ describes the characters that occur finitely often.  If we want to discuss many different types of ternary representations, then we use `or'.  That is, the notation $[\cdot,\cdot]\vee [\cdot,\cdot]\vee...\vee[\cdot,\cdot]$ is to be read as \textit{$[\cdot,\cdot]$ or $[\cdot,\cdot]$ or ... or $[\cdot,\cdot]$}. If the collection of characters occurring finitely often is empty, then we denote the corresponding type of ternary representation by $\tern{\cdot}{\emptyset}$.
\end{notation}

\begin{example}
The value $\frac{1}{2}$ has a ternary representation of $\overline{c}$.  Hence, $\frac{1}{2}$ has a type of ternary representation given by $\tern{c}{\emptyset}$.  Moreover, the value $\frac{7}{12}$ has a ternary representation given by $c\overline{rl}$, which means that $\frac{7}{12}$ has a type of ternary representation given by $\tern{lr}{c}$.
\end{example}

We note that ``the'' type of representation of a point $x\in I$ is not unique, in general.  For instance, the value $\frac{1}{3}$ has a ternary representation of type $\tern{r}{l}$ or $\tern{l}{c}$.

A thorough understanding of the ternary Cantor set is not only important for understanding many of the results on the Koch snowflake prefractal and fractal billiard. In general, Cantor sets will be ever-present and instrumental in our analysis of other fractal billiard tables.  In each example of a fractal billiard, we will clearly indicate where and how a particular Cantor set is important in analyzing a particular fractal billiard table.

\subsection{The Koch curve and Koch snowflake}
\label{subsec:TheKochCurveAndKochSnowflake}
\begin{figure}
\begin{center}
\includegraphics[scale=.5]{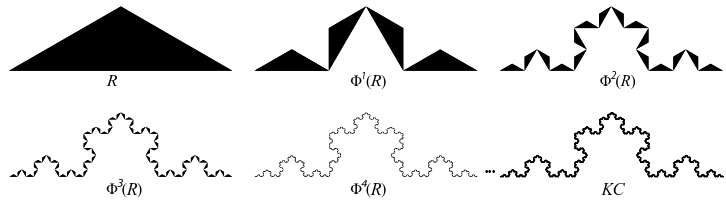}
\end{center}
\caption{The construction of the Koch curve $\kc$. Here, the self-similar set $\kc$ is viewed as a limit of the prefractal approximations $\{\Phi^m(R)\}_{m=0}^\infty$, where $R$ is the initial triangle and the map $\Phi$ is defined in terms of the IFS given by Equation (\ref{eqn:IFSForTheKochCurve}), as in Definition \ref{def:iteratedFunctionSystemAndSelfSimilarSet}. (See Theorem \ref{thm:FalconersTheorem} and the text preceding it.)}
\label{fig:ConstructionOfKochCurveFrom2DCompactSet}
\end{figure}

The Koch curve $\kc$ is constructed as shown in Figure \ref{fig:ConstructionOfKochCurveFrom2DCompactSet} and is the unique fixed point attractor of the following iterated function system on the Euclidean plane (here, $i=\sqrt{-1}$):
\begin{align}
\phi_1(\mathbf{x}) = \frac{1}{3}\mathbf{x}, \quad\quad\quad\quad\quad\quad\quad\quad\quad&
\phi_2(\mathbf{x}) = \frac{1}{3}e^{i\frac{\pi}{3}}\mathbf{x}+(\frac{1}{3},0),\label{eqn:IFSForTheKochCurve}\\
\notag\phi_3(\mathbf{x}) = \frac{1}{3}e^{-i\frac{\pi}{3}}\mathbf{x}+(\frac{2}{3},\frac{\sqrt{3}}{6}), \quad\quad\,& \phi_4(\mathbf{x}) = \frac{1}{3}\mathbf{x}+(\frac{2}{3},0).
\end{align}
\noindent Since each contraction map in the iterated function system is a similarity transformation (i.e., $\{\phi_j\}_{j=1}^4$ is a self-similar system) and $\kc = \Phi(\kc)$, we have that $\kc$ is a self-similar set; see part (ii) of Definition \ref{def:iteratedFunctionSystemAndSelfSimilarSet}.  There are additional properties of the Koch curve that are reminiscent of the Cantor set; this is more than just a coincidence and is discussed in more detail below.

\begin{figure}
\begin{center}
\includegraphics[scale=.8]{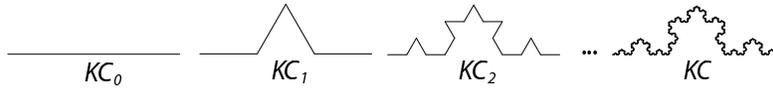}
\end{center}
\caption{One typically sees this construction of the Koch curve $\kc$ when learning about fractal sets.  Beginning with the unit interval $I$, one removes the middle third and replaces it with the two other sides of an equilateral triangle, as shown.  One then repeats this process infinitely often for every remaining interval; the resulting limiting set is $\kc$. Such a sequence $\{\kci{n}\}_{n=0}^\infty$ of approximations converges to $\kc$, because it is a subsequence of the convergent sequence of prefractal approximations $\{\Phi^m(R)\}_{m=0}^\infty$ shown in Figure \ref{fig:ConstructionOfKochCurveFrom2DCompactSet}. (Here, we are using the notion of convergence in the sense of the Hausdorff metric.)}
\label{fig:KochCurveConstruction}
\end{figure}

If we allow the iterated function system to act on the triangle $R=\{(x,y)|0\leq x\leq \frac{1}{2},0\leq y\leq \frac{\sqrt{3}}{6}x\}\cup\{(x,y)|\frac{1}{2}\leq x\leq 1,0\leq y\leq -\frac{\sqrt{3}}{6}x+\frac{\sqrt{3}}{6}\}$, as shown in Figure \ref{fig:ConstructionOfKochCurveFrom2DCompactSet}, sequential iterates of the iterated function system very quickly produce a  prefractal that is visually indiscernible from the true limiting set. But there is a more common construction that allows us to visualize the curve $\kc$ more readily, this being depicted in Figure \ref{fig:KochCurveConstruction}. The technical caveat which we are brushing under the carpet is that each polygonal approximation shown in Figure \ref{fig:KochCurveConstruction} does not contain the Koch curve $\kc$, while each approximation in the sequence shown in Figure \ref{fig:ConstructionOfKochCurveFrom2DCompactSet} does contain $\kc$.\footnote{Recall from Definition \ref{def:iteratedFunctionSystemAndSelfSimilarSet} and Theorem \ref{thm:FalconersTheorem} that for a set $F$ to be the unique fixed point attractor of an IFS, each $F_n$ must be such that $F\subseteq \Phi(F_n)$,
 so that $\Phi(F_n)=\Phi^{n+1}(F_0)$.}

\begin{notation}
For each integer $n\geq 0$, we denote by $\kci{n}$ the $n$th (inner) polygonal approximation of the Koch curve $\kc$.
\end{notation}

Intuitively, one expects the Koch curve to have finite length, since it is the limit of a sequence of polygonal approximations. On the contrary, the Koch curve $\kc$ has infinite length, which can be seen by the following calculation given in terms of the $n$th prefractal $\kci{n}$, where $\kci{n}$ is one of the polygonal approximations indicated in Figure \ref{fig:KochCurveConstruction}:
\begin{align}
\text{length of } \kci{n} &= \left(\frac{4}{3}\right)^n.
\end{align}
\noindent Then, $\lim_{n\to\infty} \left(\frac{4}{3}\right)^n = \infty$.

\begin{figure}
\begin{center}
\includegraphics[scale=.4]{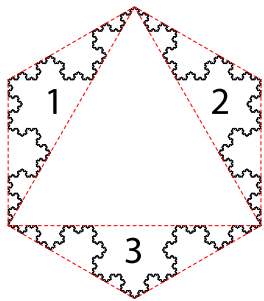}
\end{center}
\caption{The Koch snowflake is comprised of three Koch curves.  We have encapsulated each Koch curve in order to highlight how $\ks$ is the union of three abutting copies of $\kc$.}
\label{fig:3kochcurvesLabeled}
\end{figure}

The Koch snowflake $\ks$ is a fractal comprised of three abutting copies of the self-similar Koch curve; see Figure \ref{fig:3kochcurvesLabeled}.

\begin{notation}
For each integer $n\geq 0$, we denote by $\ksi{n}$ the $n$th (inner) polygonal approximation of the Koch snowflake $\ks$.
\end{notation}

As a closed (simple) curve, the Koch snowflake $\ks$ bounds a region of the plane; furthermore, the area of this region can be calculated as follows:
\begin{align}
\text{area bounded by } \ksi{n} &= 1+\sum_{i=0}^n \left(\frac{2}{3}\right)^i.
\label{eqn:AreaOfKochSnowflake}
\end{align}
\noindent Then, as $n$ increases, the right-hand side of (\ref{eqn:AreaOfKochSnowflake}) tends to a finite value. The area bounded by the Koch snowflake is thus given by $\lim_{n\to\infty} 1+\sum_{i=0}^n \left(\frac{2}{3}\right)^i = 3$, assuming the sides of $\ksi{0}$ have length one.

As we noted at the end of \S\ref{subsec:CantorSets}, Cantor sets are ever present in the context of self-similarity.  In the case of the Koch snowflake, $\ks\cap \ksi{n}$ is the union of $3\cdot 4^n$ self-similar ternary Cantor sets, each spanning a distance of $\frac{1}{3^n}$.  Such a fact will be important in determining certain sequences of what we will call \textit{compatible orbits} (see Definitions \ref{def:compatibleInitialConditions}--\ref{def:SequenceOfCompatibleOrbits}) and certain families of well-defined orbits of $\omegaks$.
\vspace{-0.25 mm}
\subsection{The $T$-fractal}
\label{subsec:TheTFractal}
\begin{figure}
\begin{center}
\includegraphics[scale=.73]{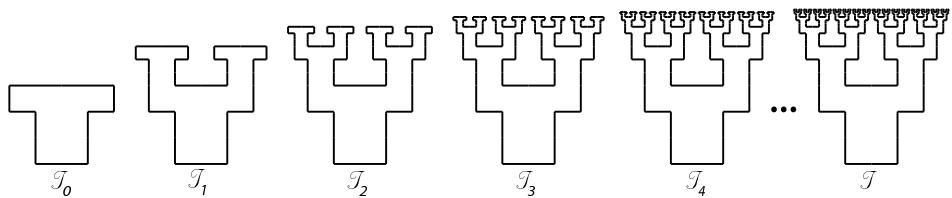}
\end{center}
\caption{The construction of the $T$-fractal $\mathscr{T}$.}
\label{fig:T-fractal}
\end{figure}

The $T$-fractal $\mathscr{T}$, discussed in \cite{AcST} in a different context, is not a self-similar set. However, $\mathscr{T}$ contains, as a proper subset, a set that is constructed in a way that is reminiscent of an iterated function system acting on a compact set so as to produce a self-similar set.\footnote{Recall from Definition \ref{def:iteratedFunctionSystemAndSelfSimilarSet} and Theorem \ref{thm:FalconersTheorem} that each prefractal approximation $F_n$ must contain the unique fixed point attractor $F$.}  As shown in Figure \ref{fig:T-fractal}, one constructs the $T$-fractal by appending scaled copies of the initial $T$ shape $\mathscr{T}_0$ to each successive approximation. Specifically, $\mathscr{T}_{n+1}$ is constructed from $\mathscr{T}_n$ by appropriately appending $2^{n+1}$ copies of $\frac{1}{2^{n+1}}\mathscr{T}_0$ to $\mathscr{T}_{n}$.\footnote{See \cite[\S{2.1}]{AcST} and \cite{LapNie6} for a more precise description of the definition of $\tfractal$.}

The overall height of $\mathscr{T}$ can be calculated and the total area bounded by $\tfractal$ can be shown to be finite, as shown in the following calculations (we assume here that the base of $\mathscr{T}_0$ is two units in length):

\begin{align}
\text{height of }\mathscr{T}_n = 3+\frac{3}{2}+\frac{3}{4}&+...+\frac{3}{2^n} = 3\sum_{i=0}^n \frac{1}{2^i}.
\end{align}

\begin{figure}
\begin{center}
\includegraphics[scale=1.2]{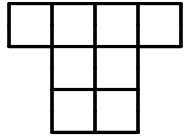}
\end{center}
\caption{$\Omega(\mathscr{T}_0)$ can be tiled by the unit square $Q$.}
\label{fig:T0TiledByUnitSquare}
\end{figure}

\noindent Then, $\lim_{n\to\infty} 3\sum_{i=0}^n \frac{1}{2^1}=6$, which is the height of $\mathscr{T}$.  Furthermore, the area bounded by $\mathscr{T}_n$ is calculated as follows. There are eight squares, each with side-length one, comprising $\mathscr{T}_0$; see Figure \ref{fig:T0TiledByUnitSquare}.  Hence, the area of $\mathscr{T}_0$ is eight square-units. Therefore,
\begin{align}
\text{area bounded by } \mathscr{T}_n = 8+ 2\cdot \frac{8}{4}+&...+2^n\cdot \frac{8}{4^n}=8\sum_{i=0}^n \frac{1}{2^i}.
\end{align}
Then, $\lim_{n\to\infty} 8\sum_{i=0}^n \frac{1}{2^i} = 16$, which is the total area bounded by $\mathscr{T}$.

There is a natural fractal subset of $\mathscr{T}$, but, for each $n\geq 0$, no point of $\mathscr{T}_n$ is in this fractal subset, which is unlike what we have seen in the case of the Koch snowflake fractal $\ks$.  In fact, the fractal subset in question is given by $\{(x,6)|x\in \mathbb{R}\}\cap \mathscr{T}$.  We note that this fractal subset is not self-similar and each point a priori fails to yield a well-defined tangent necessary for calculating the angle of reflection of a billiard ball traversing the billiard table.

\subsection{Self-similar Sierpinski carpets}
\label{subsec:aSierpinskiCarpet}
\begin{figure}
\begin{center}
\includegraphics[scale=.7]{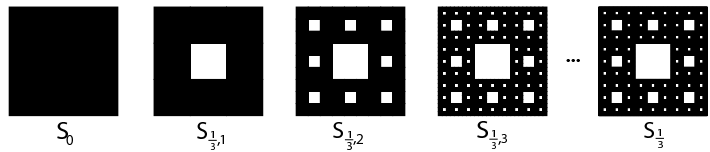}
\end{center}
\caption{The $1/3$-Sierpinski carpet is a self-similar carpet constructed in one of two ways: 1) by tremas and 2) an iterated function system (in fact, a self-similar system).  We describe here the construction of the $1/3$-Sierpinski carpet by tremas, the latter being further discussed in the main text.  Beginning with the unit square, one then removes the middle open square with side-length $\frac{1}{3}$.  From each remaining subsquare of side-length $\frac{1}{3}$, one then removes the middle open square of side-length $\frac{1}{9}$.  One continues this procedure of removing subsquares of remaining squares until there is no area left.  As one would expect, each step of the construction process can be emulated by applying the correct iterated function system, which is given in Equation (\ref{eqn:theIFSfor1-3SierpinskiCarpet}).}
\label{fig:polygonalAppxOf1-3SierpinskiCarpet}
\end{figure}

A Sierpinski carpet can be constructed by systematically removing particular open subsquares from the unit square $Q=\{(x,y)|0\leq x\leq 1,0\leq y\leq 1\}$.  Depending on how one chooses the sizes of the open subsquares to be removed, one can either construct a self-similar Sierpinski carpet or a non-self-similar Sierpinski carpet, these being defined below.  This method of construction is called \textit{construction by tremas} and is described in the caption of Figure \ref{fig:polygonalAppxOf1-3SierpinskiCarpet}, using the standard ``$1/3$-Sierpinski carpet'' as an example.  Such a construction process should be very familiar, since ``removing middle thirds'' is exactly what we did to construct the ternary Cantor set in \S\ref{subsec:CantorSets}.

As referred to in the caption of Figure \ref{fig:polygonalAppxOf1-3SierpinskiCarpet}, one may also construct the $1/3$-Sierpinski carpet by applying an appropriately defined iterated function system to the unit square $Q$.  Consider the following iterated function system, which is a self-similar system.
\begin{align}
\label{eqn:theIFSfor1-3SierpinskiCarpet}
\phi_1(\mathbf{x}) = \frac{1}{3}\mathbf{x}, & \quad\quad \phi_2(\mathbf{x}) = \frac{1}{3}\mathbf{x}  +\left(0,\frac{1}{3}\right), \\
\notag\phi_3(\mathbf{x}) = \frac{1}{3}\mathbf{x} +\left(0,\frac{2}{3}\right), & \quad\quad \phi_4(\mathbf{x}) = \frac{1}{3}\mathbf{x}  +\left(\frac{1}{3},0\right),\\
\notag\phi_5(\mathbf{x}) = \frac{1}{3}\mathbf{x} +\left(\frac{1}{3},\frac{2}{3}\right), & \quad\quad\phi_6(\mathbf{x}) = \frac{1}{3}\mathbf{x}  +\left(\frac{2}{3},0\right),\\
\notag\phi_7(\mathbf{x}) = \frac{1}{3}\mathbf{x} +\left(\frac{2}{3},\frac{1}{3}\right), &\quad\quad \phi_8(\mathbf{x}) = \frac{1}{3}\mathbf{x}  +\left(\frac{2}{3},\frac{2}{3}\right).
\end{align}
\noindent Then, denoting the $1/3$-Sierpinski carpet by $\Sa{3}$, we have that $\lim_{n\to\infty} \Phi^n(Q) = \Sa{3}$.  Since each contraction in the iterated function system is a similarity contraction and $\Sa{3} = \Phi(\Sa{3})$, it follows that $\Sa{3}$  is a self-similar set.

We discuss here the relevant results and material from \cite{Du-CaTy}.  For our purposes, the first level approximation of a Sierpinski carpet $S_{\mathbf{a}}$ will always be the unit square $Q$ and denoted by $S_0$.   Since every (self-similar and non-self-similar) Sierpinski carpet has the same zeroth level approximation and zero is never a scaling ratio, such notation will never cause any confusion.

What follows is a general description on how to construct a Sierpinski carpet by removing appropriately sized middle open squares.  Consider the unit square $Q=S_0$.  Let $a_0 = 2k_0+1$ for some $k_0\in \N$.  Partition $S_0$ into $a_0$ squares of side-length $a_0^{-1}$.  Next, remove the middle open subsquare.  Let $a_1 =2k_1 + 1$ for some $k_1\in \N$. Each subsquare may then be partitioned into $a_1^2$ many squares with side-length $(a_0\cdot a_1)^{-1}$.  We then remove each middle open subsquare of side-length $(a_0\cdot a_1)^{-1}$; see Figure \ref{fig:polygonalAppxOf1-3SierpinskiCarpet}.  Continuing this process, let $a_{n-1} = 2k_{n-1}+1$ where $k_{n-1}\in \N$ and let $a_n = 2k_n+1$ for some $k_n\in N$.  Then we partition a subsquare of side-length $(a_0\cdot a_1\cdots a_{n-1})^{-1}$ into $a_{n}^2$ many squares.  We then remove the middle open square from each subsquare in the partition.  Continuing in this manner ad infinitum, one constructs a Sierpinski carpet denoted by $\Sa{a}$, where $\mbfa = \{a_i^{-1}\}_{i=0}^\infty$.


\begin{definition}[A self-similar Sierpinski carpet]
\label{def:ASelfSimilarSierpinskiCarpet}
If $\mathbf{a}=\{a_i^{-1}\}_{i=0}^\infty$, with $a_i = 2k_i+1$ and $k_i \in \N$, is a periodic sequence of rational values, then the Sierpinski carpet $S_\mbfa$ is called a \textit{self-similar Sierpinski carpet}.
\end{definition}

We have described the construction of a self-similar Sierpinski carpet $S_{\mbfa}$ in terms of the removal of particular open squares.  As the name would suggest, there exists a suitably defined iterated function system $\{\phi_i\}_{i=1}^k$ such that $S_\mbfa = \Phi(S_\mbfa)$.  Viewing $S_\mathbf{a}$ as the unique fixed point attractor of an appropriately defined iterated function system will be useful in stating some of the results in the subsequent sections.  More precisely, $S_{\mbfa}$ is viewed as the self-similar set associated with a self-similar system, as in part (ii) of Definition \ref{def:iteratedFunctionSystemAndSelfSimilarSet}.

While we do not discuss any results concerning non-self-similar Sierpinski carpet billiards in this paper, we provide the definition for completeness.

\begin{definition}[A non-self-similar Sierpinski carpet]
If $\mathbf{a}=\{a_i^{-1}\}_{i=0}^\infty$, with $a_i = 2k_i+1$ and $k_i \in \N$, is an aperiodic sequence of rational values, then the Sierpinski carpet $S_\mbfa$ is called a \textit{non-self-similar Sierpinski carpet}.
\end{definition}

\begin{definition}[A cell of $S_{\mathbf{a},n}$]
\label{def:ACellOfSai}
Let $a_0=2k_0+1$, $k_0\in \mathbb{N}$.  Consider a partition of the unit square $Q=S_0$ into $a_0^2$ many squares of side-length $a_0^{-1}$.  A subsquare of the partition is called a \textit{cell of $S_0$} and is denoted by $\celli{0}{a_0}$.  Furthermore, let $S_\mathbf{a}$ be a Sierpinski carpet.  Consider a partition of the prefractal approximation $S_{\mathbf{a},n}$ into subsquares with side-length $(a_0\cdot a_1\cdots a_n)^{-1}$.  A subsquare of the partition of $S_{\mathbf{a},n}$ is called a \textit{cell of $S_{\mathbf{a},n}$} and is denoted by $\celli{n}{a_0a_1\cdots a_n}$ and has side-length $(a_0\cdot a_1\cdots a_n)^{-1}$.
\end{definition}

\begin{definition}[Peripheral square]
In accordance with the convention adopted in \cite{Du-CaTy}, the boundary of an open square removed in the construction of $S_\mbfa$ is called a \textit{peripheral square} of $S_\mbfa$.  Furthermore, by convention, the unit square $Q=S_0$ is not a peripheral square.
\end{definition}

\begin{definition}[Nontrivial line segment of $S_\mbfa$]
\label{def:nontrivialLineSegmentOfSa}
A \textit{nontrivial line segment of $S_\mbfa$} is a (straight-line) segment of the plane contained in $S_\mbfa$ and which has nonzero length.
\end{definition}

Unless otherwise indicated, in what follows, we assume that $S_\mbfa$ is a self-similar Sierpinski carpet with a single scaling ratio $a$; that is, $\mathbf{a}=\{a^{-1}\}_{i=0}^\infty$, where $a = 2k+1$ for some fixed $k\in \N$.  In addition, when $\mbfa = \{a^{-1}\}_{i=0}^\infty$, $S_{\mbfa}$ is denoted by $S_a$.

We next state the following theorem, due to Durand-Cartagena and Tyson in \cite{Du-CaTy} and which will be very useful to us in this context  (see \S\ref{subsec:AprefractalSelfSimilarSierpinskiCarpetBilliard} and \S\ref{subsec:aSelf-SimilarSierpinskiCarpetBilliard}).

\begin{theorem}[{\cite[Thm. 4.1]{Du-CaTy}}]
\label{thm:AsetBset}
Let $S_a$ be a self-similar Sierpinski carpet.  Then the set of slopes $\slopesa{a}$ of nontrivial line segments of $S_a$ is the union of the following two sets\emph{:}
\begin{align}
A &= \left\{\frac{p}{q}: p+q\leq a,\,\,0\leq p<q\leq a-1,\,\, p,q\in \mathbb{N}\cup\{0\},\,\,p+q \text{ is odd}\right\}, \label{eqn:Aset}\\
B &= \left\{\frac{p}{q}: p+q\leq a-1,\,\,0\leq p\leq q\leq a-2,\,\, p,q\in \mathbb{N},\,\,p,q \text{ are odd}\right\}.\label{eqn:Bset}
\end{align}

Moreover, if $\alpha\in A$, then each nontrivial line segment in $S_a$ with slope $\alpha$ touches vertices of peripheral squares, while if $\alpha\in B$, then each nontrivial line segment in $S_a$ with slope $\alpha$ is disjoint from all peripheral squares.
\end{theorem}

\begin{notation}
\label{not:AaBaAbBb}
Let $a,b$ be odd positive integers such that $3\leq b\leq a$ and let $\slopesa{a}$ and $\slopesa{b}$ be the set of slopes of nontrivial line segments of $S_a$ and $S_b$, respectively.  We denote by $A_a$ (resp., $A_b$) the subset $A\subseteq \slopesa{a}$ (resp., $A\subseteq\slopesa{b}$) given in Equation (\ref{eqn:Aset}) of Theorem \ref{thm:AsetBset}.  Similarly, we denote by $B_a$ (resp., $B_b$) the subset $B\subseteq \slopesa{a}$ (resp., $\slopesa{b}$) given in Equation (\ref{eqn:Bset}) of Theorem \ref{thm:AsetBset}.\footnote{In the case of $A_b$ (resp., $B_b$), $a$ should of course be replaced by $b$ in Equation (\ref{eqn:Aset}) (resp., Equation (\ref{eqn:Bset})).}
\end{notation}

If $S_a$ and $S_b$ are self-similar Sierpinski carpets with $b\leq a$, then it is clear that $\slopesa{b}\subseteq\slopesa{a}$.  Moreover, in this case, we also have that $A_b\subseteq A_a$ and $B_b\subseteq B_a$.

\begin{remark}
We note that if $\alpha$ is the slope of a nontrivial line segment in $S_a$, then so is $-\alpha$, $\frac{1}{\alpha}$ and $-\frac{1}{\alpha}$ by symmetry of the carpet.  However, we restrict our attention in this paper to the slopes described in the above result of \cite{Du-CaTy}.
\end{remark}

\section{Prefractal (rational) billiards}
\label{sec:PrefractalRationalBilliards}

In the previous sections, we surveyed basic facts and results from mathematical billiards and fractal geometry, with most of our attention being focused on the subject of rational billiards and sets exhibiting self-similarity.  We also discussed the importance of examining the dynamically equivalent geodesic flow on an associated translation surface.  In this section, we will examine examples from particular classes of prefractal (rational) billiards.  We are interested in tables that can be tiled by a single polygon which can also tile the (Euclidean) plane.  The main examples we will discuss are the Koch snowflake prefractal billiard table, the $T$-fractal prefractal billiard table and a self-similar Sierpinski carpet prefractal billiard table.  Each example of a prefractal billiard table constitutes a rational billiard table, but is an element of a sequence of rational billiard tables approximating a fractal billiard table with radically different qualities when compared to the others.  That is, the Koch snowflake has an everywhere nondifferentiable boundary; the $T$-fractal billiard table is certainly a fractal billiard table, since its boundary $\mathscr{T}$ contains a fractal set, but the portion of the boundary that is nondifferentiable has Lebesgue measure zero; a Sierpinski carpet billiard table can possibly have no area, yet yield billiard orbits of finite length.

\subsection{A general structure}
We restrict our attention to billiard tables with fractal boundary $F$, where $F$ can be approximated by a suitably chosen sequence of rational polygons $\{F_n\}_{n=0}^\infty$.  More specifically, we are interested in a fractal billiard table $\Omega(F)$ with the property that, for every $n\geq 0$, $\Omega(F_n)$ can be tiled by a single polygon $D_n$, where $D_n =  c_n D_0$.  Here, $0<c_n\leq 1$ is a suitably chosen scaling ratio and $D_0$ is a polygon that tiles both the (Euclidean) plane as well as the rational billiard $\omegafraci{0}$.\footnote{In the case of certain prefractal approximations, $\omegafraci{0}$ is exactly $D_0$.  In general, however, $D_0$ does not always equal $\omegafraci{0}$, but certainly tiles $\omegafraci{0}$. An example of this situation is $\omegati{0}$; see \S\ref{subsec:TheTFractalPrefractalBilliard}.  Such a billiard is tiled by the unit square, which is the associated polygon $D_0$.}

The focus in this subsection is on developing a general framework for discussing billiards on prefractal approximations.  If $\omegafraci{n}$ and $\omegafraci{n+1}$ are two prefractal billiard tables approximating a given fractal billiard table $\omegafrac$, then we want to have a systematic way of determining how and if two orbits $\ofraci{n}$ and $\ofraci{n+1}$ of $\omegafraci{n}$ and $\omegafraci{n+1}$ are related.

\begin{notation}
We will primarily measure angles relative to a fixed coordinate system, with the origin being fixed at a corner of a prefractal approximation $F_0$.  However, we will sometimes measure an angle relative to a side of $F_n$ on which a billiard ball lies.  In such situations, we will write the angle as $\varpi(\theta)$ in order to indicate that the inward pointing direction is $\theta$, measured relative to the side on which the vector is based.
\end{notation}

\begin{figure}
\begin{center}
\includegraphics[scale=.6]{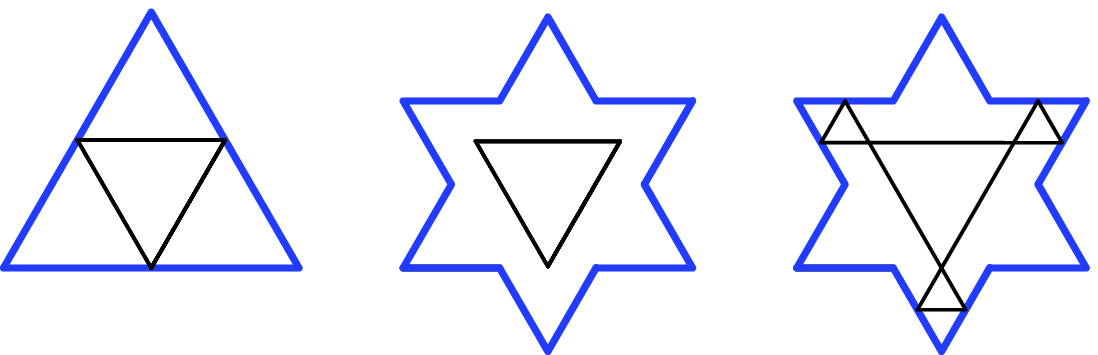}
\end{center}
\caption{In the first image, we have the orbit $\ofracixang{0}{\overline{c}}{\frac{\pi}{3}}$ of $\omegaksi{0}$. In the second image, we see that the orbit $\ofracixang{0}{\overline{c}}{\frac{\pi}{3}}$, when embedded in $\omegaksi{1}$, is not an orbit of $\omegaksi{1}$. In the third image, the given orbit of $\omegaksi{1}$ intersects sides of $\omegaksi{1}$ and appears to be related to $\ofracixang{0}{\overline{c}}{\frac{\pi}{3}}$ in some way.}
\label{fig:notApwFagnanoOrbit}
\end{figure}

To motivate our general discussion, consider the orbit $\ofraciang{0}{\frac{\pi}{3}}$ of $\omegaksi{0}$, where $\xoo = \overline{c}\in I$; see the first image in Figure \ref{fig:notApwFagnanoOrbit} (and recall our earlier discussion in \S\ref{subsec:CantorSets}).  The same orbit, viewed as a continuous curve embedded in $\omegaksi{1}$, does not constitute an orbit of $\omegaksi{1}$; see the second image in Figure \ref{fig:notApwFagnanoOrbit}.  Consider the orbit $\ofraciang{1}{\frac{\pi}{3}}$ shown in the third image in Figure \ref{fig:notApwFagnanoOrbit}.  Such an orbit does intersect the boundary of $\omegaksi{1}$ and appears to be related to $\ofracixang{0}{\overline{c}}{\frac{\pi}{3}}$, but in what way we have not yet explicitly said.  Initially, we notice that, as a continuous curve embedded in $\omegaksi{1}$, the orbit $\ofracixang{0}{\overline{c}}{\frac{\pi}{3}}$ is a subset of $\ofraciang{1}{\frac{\pi}{3}}$.  Being eager to establish a proper notion of ``related'', we may be inclined to declare that two orbits are related if one is a subset of the other, when viewed as continuous curves in the plane.  Unfortunately, we quickly see that such a definition is highly restrictive.  A more general observation is that $\xio{1}$ and $\xoo$ are collinear in the direction of $\frac{\pi}{3}$, without any portion of $\ksi{1}$ intersecting the segment $\overline{\xio{1}\xoo}$.  We then say that $(\xio{0},\frac{\pi}{3})$ and $(\xio{1},\frac{\pi}{3})$ are \textit{compatible initial conditions}.  We state the formal definition as follows.


\begin{definition}[Compatible initial conditions]
\label{def:compatibleInitialConditions}
Without loss of generality, suppose that $n$ and $m$ are nonnegative integers such that $n > m$. Let $(\xio{n},\theta_n^0)\in (\Omega(F_n)\times S^1)/\sim$ and $(\xio{m},\theta_m^0)\in (\Omega(F_m)\times S^1)/\sim$ be two initial conditions of the orbits $\ofraci{n}$ and $\ofraci{m}$, respectively, where we are assuming that $\theta_n^0$ and $\theta_m^0$ are both inward pointing.  If $\theta_n^0 = \theta_m^0$ and if $\xio{n}$ and $\xio{m}$ lie on a segment determined from $\theta_n^0$ (or $\theta_m^0$) that intersects $\Omega(F_n)$ only at $\xio{n}$, then we say that $(\xio{n},\theta_n^0)$ and $(\xio{m},\theta_m^0)$ are \textit{compatible initial conditions}.
\end{definition}

\begin{remark}
When two initial conditions $(\xio{n},\theta_n^0)$ and $(\xio{m},\theta_m^0)$ are compatible, then we simply write each as $(\xio{n},\theta^0)$ and $(\xio{m},\theta^0)$. If two orbits $\ofraci{m}$ and $\ofraci{n}$ have compatible initial conditions, then we say such orbits are \textit{compatible}.
\end{remark}

Depending on the nature of $\Omega(F)$, not every orbit must pass through the region of $\omegafraci{n}$ corresponding to the interior of $\omegafraci{0}$, let alone pass through the interior of $\omegafraci{m}$, for any $m<n$.  Because of this, it may be the case that an initial condition $(\xio{n},\theta^0)$ is not compatible with $(\xio{m},\theta^0)$, for any $m<n$.  As such, in Definitions \ref{def:sequenceOfCompatibleInitialConditions} and \ref{def:SequenceOfCompatibleOrbits}, we consider sequences beginning at $i=N$, for some $N\geq 0$.

\begin{definition}[Sequence of compatible initial conditions]
\label{def:sequenceOfCompatibleInitialConditions}
Let $\{(\xio{i},\theta_i^0)\}_{i=N}^\infty$ be a sequence of initial conditions, for some integer $N\geq 0$.  We say that this sequence is a \textit{sequence of compatible initial conditions} if for every $m\geq N$ and for every $n> m$, we have that $(\xio{n},\theta_n^0)$ and $(\xio{m},\theta_m^0)$ are compatible initial conditions.  In such a case, we then write the sequence as $\{(\xio{i},\theta^0)\}_{i=N}^\infty$.
\end{definition}

\begin{definition}[Sequence of compatible orbits]
\label{def:SequenceOfCompatibleOrbits}
Consider a sequence of compatible initial conditions $\{(\xio{n},\theta^0)\}_{n=N}^\infty$.  Then the corresponding sequence of orbits $\compseqi{N}$ is called \textit{a sequence of compatible orbits}.
\end{definition}

If $\ofracixang{m}{\xio{m}}{\theta_m^0}$ is an orbit of $\omegafraci{m}$, then $\ofracixang{m}{\xio{m}}{\theta_m^0}$ is a member of a sequence of compatible orbits $\compseqiang{N}{\theta^0}$ for some $N\geq 0$. It is clear from the definition of a sequence of compatible orbits that such a sequence is uniquely determined by the first orbit $\ofraciang{N}{\theta^0}$.  Since the initial condition of an orbit determines the orbit, we can say without any ambiguity that a sequence of compatible orbits is determined by an initial condition $(\xio{N},\theta^0)$.

\begin{definition}[A sequence of compatible $\mathcal{P}$ orbits]
Let $\mathcal{P}$ be a property (resp., $\mathcal{P}_1,...,\mathcal{P}_j$ a list of properties).  If every orbit in a sequence of compatible orbits has the property $\mathcal{P}$ (resp., a list of properties $\mathcal{P}_1,...,\mathcal{P}_j$), then we call such a sequence \textit{a sequence of compatible $\mathcal{P}$ \emph{(}resp., $\mathcal{P}_1,...,\mathcal{P}_j$\emph{)}  orbits}.
\end{definition}


The following theorem can be deduced from Theorem 3 of Gutkin's paper \cite{Gut2}; see \cite[\S{3.2}]{LapNie3}.

\begin{theorem}
\label{thm:topologicalDichotomyForFn}
Consider a prefractal rational billiard $\omegafraci{n}$. If $\Omega(\fraci{n})$ is tiled by a rational polygon $D_n$ such that $D_n$ tiles the Euclidean plane, then, for a fixed direction $\theta_n^0$, every orbit $\ofraciang{n}{\theta_n^0}$ of $\omegafraci{n}$ is closed or every orbit $\ofraciang{n}{\theta_n^0}$ is dense in $\omegafraci{n}$,\footnote{Recall that these notions were introduced towards the beginning of \S\ref{sec:RationalBilliards}.} regardless of the initial basepoint $\xio{n}$.
\end{theorem}

\begin{remark}
When $\omegafraci{n}$ is tiled by $D_n$, where $D_n$ is a rational polygon tiling the plane, then $\omegafraci{n}$ is more generally referred to as an \textit{almost integrable billiard}, this being the language used in \cite{Gut2}.
\end{remark}

The following is a generalization to this broader setting of Corollary 16 from \cite{LapNie3}.  It is established in the same manner.

\begin{theorem}
\label{thm:generalTopologicalDichotomyForSequencesOfCompatibleOrbits}
Let $\Omega(F)$ be a fractal billiard table approximated by a suitable sequence of rational polygonal billiard tables $\{\Omega(F_n)\}_{n=0}^\infty$.  If there exists a polygon $D_0$ that tiles the plane and such that for every $n\geq 0$ there exists $0<c_n\leq 1$ with $D_n:=c_n D_0$ tiling $\omegafraci{n}$, then any sequence of compatible orbits is either entirely comprised of closed orbits or entirely comprised of orbits that are dense in their respective billiard tables.
\end{theorem}


\subsection{The prefractal Koch snowflake billiard}
\label{subsec:ThePrefractalKochSnowflakeBilliard}

The billiard $\omegaksi{n}$ can be tiled by equilateral triangles.  Specifically, if $\Delta$ is the equilateral triangle with sides having unit length, then $\omegaksi{n}$ is tiled by $\frac{1}{3^n}\Delta$, for every $n\geq 0$.  Moreover, as is well known, $\Delta = \ksi{0}$ tiles the plane.  Therefore, Theorems \ref{thm:topologicalDichotomyForFn} and \ref{thm:generalTopologicalDichotomyForSequencesOfCompatibleOrbits} hold for the prefractal billiard $\omegaksi{n}$.

Our goal for this subsection and \S\ref{subsubsec:TheCorrespondingPrefractalTranslationSurface} is to survey some of the main results of \cite{LapNie1, LapNie2, LapNie3}.  We will focus on pertinent examples that will motivate a richer discussion in \S\ref{subsec:TheKochSnowflakeFractalBilliard}. Initially, we focus on properties of orbits with an initial direction of $\frac{\pi}{3}$ and $\frac{\pi}{6}$.\footnote{Equivalently, we could focus on orbits with an initial direction of $\frac{\pi}{3}$ and $\frac{\pi}{2}$, since $\frac{\pi}{2}$ is the rotation of $\frac{\pi}{6}$ through the angle $\frac{\pi}{3}$, the angle $\frac{\pi}{3}$ being an angle that determines an axis of symmetry of $\ksi{n}$, for $n\geq 0$.}

If $\ofraciang{0}{\frac{\pi}{3}}$ is an orbit of $\omegaksi{0}$, so long as $\xio{0}$ is not a corner of $\ksi{0}$, the orbit will be periodic, as expected.  However, depending on the nature of the ternary representation of $\xio{0}$, the compatible orbit $\ofraciang{1}{\frac{\pi}{3}}$ may be singular in $\omegaksi{1}$.\footnote{Recall that the notion of type of a ternary representation was introduced in Notation \ref{not:typeOfTernaryRepresentation} of \S\ref{subsec:CantorSets}.}

\begin{theorem}[{\cite{LapNie3}}]
Let $\xio{0}\in I\subseteq\ksi{0}$.  If $\xoo$ has a ternary representation of type $\tern{l}{cr}\vee \tern{r}{lc}$, then there exists $N\geq 0$ such that the compatible orbit $\ofraciang{N}{\frac{\pi}{3}}$ will be singular in $\omegaksi{N}$.  Moreover, for every $n\geq N$, $\ofraciang{n}{\frac{\pi}{3}}$ will also be singular in $\omegaksi{n}$.
\end{theorem}

\begin{theorem}[{\cite{LapNie3}}]
If $\xio{0}$ has a ternary representation of the form $\tern{c}{lr}\vee\tern{lc}{r}\vee\tern{cr}{l}\vee\tern{lcr}{\emptyset}\vee\tern{lr}{c}$, then the sequence of compatible orbits given by $\compseqiang{0}{\frac{\pi}{3}}$ is a sequence of compatible periodic orbits.
\end{theorem}


\begin{theorem}
The length and period of an orbit $\ofraciang{m}{\frac{\pi}{3}}\in \compseqang{\frac{\pi}{3}}$ is dictated by the ternary representation of $\xio{0}$. \emph{(}See \cite{LapNie2} for the corresponding specific formulas.\emph{)}
\end{theorem}

\begin{remark}
See \S{4.4} of \cite{LapNie2} for a precise statement and proof of this result, as well as for additional properties of orbits with an initial direction of $\frac{\pi}{3}$.
\end{remark}

\begin{figure}
\begin{center}
\includegraphics[scale=0.25]{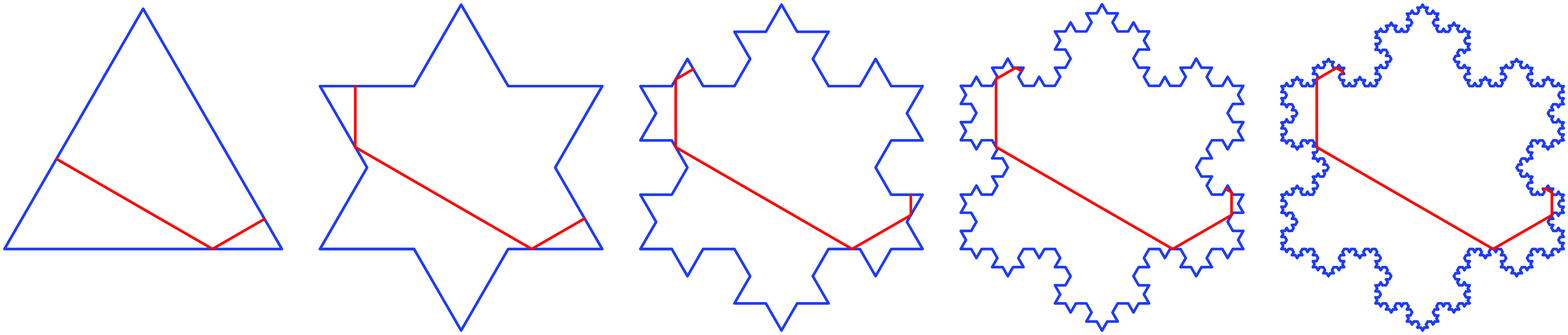}
\end{center}
\caption{An example of a hook orbit.  The same initial condition is used in each prefractal billiard.}
\label{fig:hybridpi6lv5}
\end{figure}

\begin{example}[A sequence of compatible hook orbits]
\label{exa:ASequenceOfCompatibleHookOrbits}
Let $\xoo\in I$ have a ternary representation given by $\overline{rl}$. Such a representation indicates that, in each prefractal approximation $\ksi{n}$, $\xoo$ is an element of an open, connected neighborhood contained in $\ksi{n}$. The point $\xoo$ corresponds to the value $3/4\in I$.  If we consider an orbit of $\omegaksi{0}$ with an initial direction of $\frac{\pi}{6}$, the ternary representation of the basepoints at which the billiard ball path forms right angles with the sides of $\omegaksi{0}$ is of the type $\tern{c}{lr}$. This is a degenerate periodic hybrid orbit, meaning that it doubles back on itself, and the next orbit in the sequence of compatible periodic hybrid orbits has the initial condition $(\xio{1},\frac{\pi}{6}) = (\xoo,\frac{\pi}{6})$.  Since the ternary representation of the basepoint of $f_0(\xoo,\frac{\pi}{6})$ is $r\overline{c}$ and $\theta_0^0=\theta_1^0=\frac{\pi}{6}$, it follows that the basepoint of $f_1(\xio{1},\frac{\pi}{6})$ is a point which, for every prefractal approximation $\ksi{n}$, is an element of an open, connected neighborhood contained in $\ksi{n}$.  Then the basepoint of $f^2_1(\xio{1},\frac{\pi}{6})$ (where $f_1^2$ denotes the second iterate of the billiard map $f_1$) has a ternary representation of type $\tern{c}{lr}$.  This same pattern is repeated for every subsequent orbit in the sequence of compatible orbits.  It follows that the resulting sequence of compatible orbits forms a sequence of orbits that is converging to a set which is well defined. That is, such a set will be some path in the fractal billiard table $\omegaks$ with finite length which is effectively determined by the law of reflection in each prefractal approximation of $\omegaks$.

Such orbits are introduced in \cite{LapNie3} and referred to as \textit{hook orbits}, because they appear to be ``hooking'' into the Koch snowflake; see Figure \ref{fig:hybridpi6lv5}.
\end{example}

The hook orbits of Example \ref{exa:ASequenceOfCompatibleHookOrbits} are special cases of a general class of orbits called hybrid orbits, which were introduced, as well as studied, in \cite{LapNie3}.

\begin{definition}[Hybrid orbit]
\label{def:hybridOrbit}
Let $\ofraci{n}$ be an orbit of $\omegaksi{n}$.  If all but at most two basepoints  $\xii{n}{k_n}\in \fprintksi{n}$ have ternary representations (determined with respect to the side $s_{n,\nu}$ on which each point resides) of type $\tern{c}{lr}\vee\tern{cl}{r}\vee\tern{cr}{l}\vee\tern{lcr}{\emptyset}\vee\tern{lr}{\emptyset}$, then we call $\ofraciang{n}{\theta_n^0}$ a \textit{hybrid orbit} of $\omegaksi{n}$.
\end{definition}

\begin{figure}
\begin{center}
\includegraphics[scale=.25]{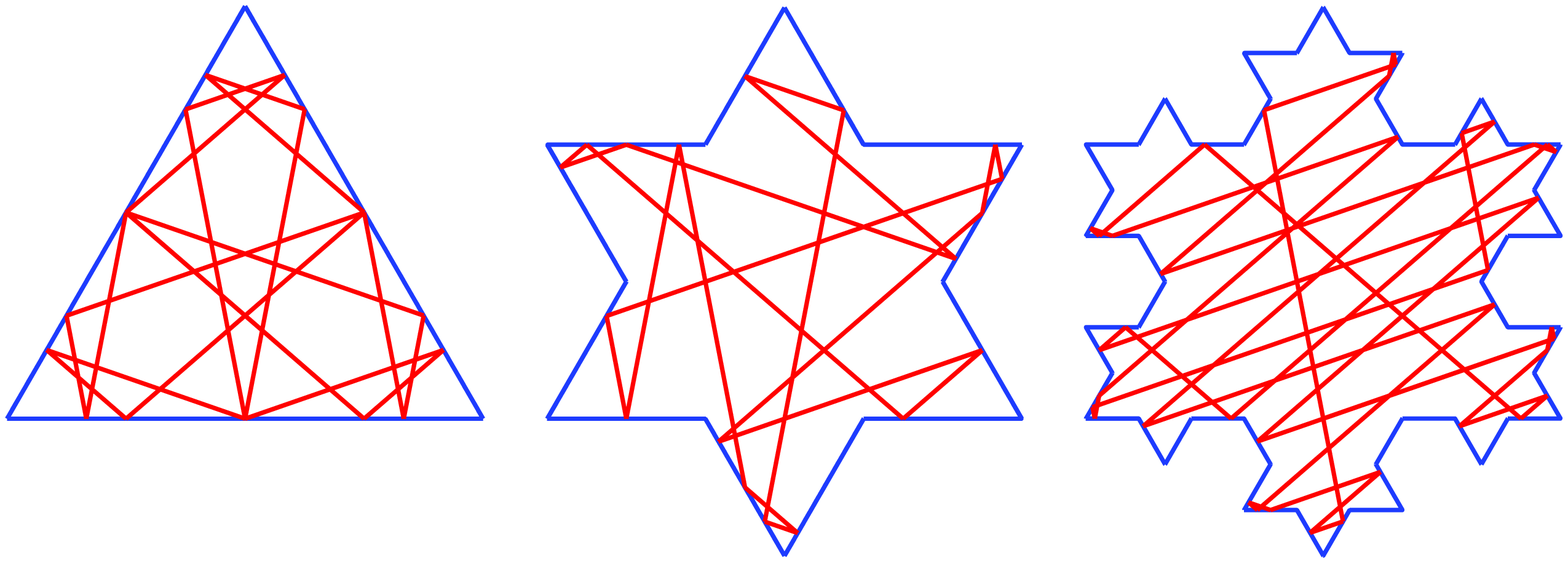}
\end{center}
\caption{Three examples of periodic hybrid orbits.  These are the first three elements of the sequence of compatible periodic hybrid orbits described in Example \ref{exa:ASequenceOfCompatiblePeriodicHybridOrbits}.  In order to understand exactly what is discussed in the paragraph immediately following Definition  \ref{def:hybridOrbit}, compare and contrast the hybrid orbits shown here with the hybrid orbits shown in Figures \ref{fig:notApwFagnanoOrbit}, \ref{fig:hybridpi6lv5} and \ref{fig:CompatibleCantorOrbit7-12}.  Certain segments of the hybrid orbits shown here remain intact and become subsets of subsequent compatible periodic hybrid orbits, yet the orbits are visually different from one another.}
\label{fig:hybrid012}
\end{figure}

A hybrid orbit is so named for the fact that it \textit{may} have qualities reminiscent of an orbit $\ofraci{n+1}$ that is identical to the compatible orbit $\ofraci{n}$ and an orbit $\ofracixang{n+1}{y_{n+1}^0}{\gamma^0_{n+1}}$ that is visually different from the compatible orbit $\ofracixang{n}{y_n^0}{\gamma^0_n}$; see Figure \ref{fig:hybrid012} and its caption.

\begin{definition}[A $\mathscr{P}$ hybrid orbit]
If $\ofraci{n}$ is a hybrid orbit with property $\mathscr{P}$, then we say that it is a $\mathscr{P}$ hybrid orbit.
\end{definition}




\begin{proposition}
\label{prop:denseOrbitIsADenseHybridOrbit}
If $\ofraci{n}$ is a dense orbit of $\omegaksi{n}$, then $\ofraci{n}$ is a dense hybrid orbit.
\end{proposition}

Applying the results in Theorem \ref{thm:generalTopologicalDichotomyForSequencesOfCompatibleOrbits} and Proposition \ref{prop:denseOrbitIsADenseHybridOrbit}, we state the following result.

\begin{theorem}[A topological dichotomy for sequences of compatible orbits of prefractal billiard tables, \cite{LapNie3}]
\label{thm:ATopologicalDichotomy}
Let $\compseqi{N}$ be a sequence of compatible orbits.  Then  we have that $\compseqi{N}$ is either entirely comprised of closed orbits or is entirely comprised of dense hybrid orbits.\footnote{Recall that the notions of ``closed orbit'' and ``dense orbit'' were defined in \S\ref{sec:RationalBilliards}, just before Definition \ref{def:footprintOfAnOrbit}.}
\end{theorem}

\begin{theorem}[{\cite{LapNie3}}]
\label{thm:hybridOrbitOfKS0ImpliesHybridOrbitOfKSn}
If $\ofraci{0}$ is a periodic hybrid orbit of $\omegaksi{0}$ with no basepoints corresponding to ternary points \emph{(}i.e., points having ternary representations of the types $\tern{l}{cr}\vee\tern{r}{lc}$\emph{)}, then for every $n\geq 0$, the compatible orbit $\ofraci{n}$ is a periodic hybrid orbit of $\omegaksi{n}$.
\end{theorem}


In order to fully understand the following result, we define what it means for a vector to be \textit{rational with respect to a basis} $\{u_1,u_2\}$ of $\mathbb{R}^2$.  If $z=mu_1+nu_2$, for some $m,n\in \mathbb{Z}$, then we say that $z$ is \textit{rational with respect to the basis} $\{u_1,u_2\}$.  Otherwise, we say that $z$ is \textit{irrational with respect to} $\{u_1,u_2\}$.

\begin{theorem}[A sequence of compatible periodic hybrid orbits, \cite{LapNie3}]
Let $\xoo\in I$ and consider a vector $(a,b)$ that is rational with respect to the basis $\{u_1,u_2\}:=\{(1,0),(1/2,\sqrt{3}/2)\}$.  Then, we have the following\emph{:}

\begin{enumerate}
\item{If $a$ and $b$ are both positive integers with $b$ being odd, $\xoo = \frac{r}{4^s}$, for some $r,s\in\mathbb{N}$ with $s\geq 1$, $1\leq r<4^s$ being odd  and $\theta^0 := \arctan{\frac{b\sqrt{3}}{2a+b}}$, then the sequence of compatible closed orbits $\compseq$ is a sequence of compatible periodic hybrid orbits.}

\item{If $a=1/2$, $b$ is a positive odd integer, $\xoo = \frac{r}{2^s}$, for some $r,s\in\mathbb{N}$ with $s\geq 1$, $1\leq r<2^s$ being odd and $\theta^0 := \arctan{\frac{b\sqrt{3}}{2a+b}}$, then the sequence of compatible closed orbits $\compseq$ is a sequence of compatible periodic hybrid orbits.}
\end{enumerate}
\label{thm:bodd}
\end{theorem}

\begin{remark}
We want to emphasize that the angle $\theta^0$ in Part (1) and Part (2) of Theorem \ref{thm:bodd} is not necessarily $\frac{\pi}{3}$, $\frac{\pi}{2}$ or $\frac{\pi}{6}$, but can assume countably infinitely many values.
\end{remark}

\begin{example}[A sequence of compatible periodic hybrid orbits]
\label{exa:ASequenceOfCompatiblePeriodicHybridOrbits}
In Figure \ref{fig:hybrid012}, three periodic hybrid orbits are displayed.  These three orbits constitute the first three terms in a sequence of compatible periodic hybrid orbits.\footnote{By Theorem \ref{thm:bodd}(2), the angle $\theta^0_0$ determined by the initial segment of the orbit and the initial basepoint $\xoo=\frac{1}{2}=\overline{c}$ both guarantee that the sequence of compatible orbits $\compseqi{0}$ is a sequence of compatible periodic hybrid orbits.}  If we choose $\xoo=\overline{c}\in I$ and $\theta_0^0$ to be an angle such that $\xoo$ connects with the midpoint of the lower one-third interval on the side of $\omegaksi{0}$, we can see that $\ofraci{0}$ is a periodic hybrid orbit.  More importantly, there are elements of the footprint $\fprintksi{0}$ with ternary representations of type $\tern{lr}{c}$.  This observation is key for constructing what we call nontrivial paths of $\omegaks$, a topic which is discussed in more detail in \S\ref{subsec:TheKochSnowflakeFractalBilliard}.
\end{example}

Given a nonnegative integer $N$, we say that a sequence of compatible orbits $\compseqi{N}$ is a \textit{constant sequence of compatible orbits} if the path traversed by $\ofraci{n+1}$ is identical to the path traversed by $\ofraci{n}$, for every $n\geq N$.  Furthermore, we say that a sequence of compatible orbits $\compseqi{0}$ is \textit{eventually constant} if there exists a nonnegative integer $N$ such that $\compseqi{N}$ is constant, in the above sense.

\begin{theorem}[A constant sequence of compatible periodic hybrid orbits, \newline \cite{LapNie3}]
\label{thm:SufficientConditionForCantorOrbit}
Let $\ofraci{0}$ be an orbit of $\omegaksi{0}$ such that every $\xii{0}{k_0}\in\fprintksi{0}$ has a ternary representation of type $\tern{lr}{c}$.  Then $\compseqi{0}$ is a sequence of compatible periodic hybrid orbits. Moreover, there exists $N\geq 0$ such that $\compseqi{N}$ is a constant sequence of compatible periodic hybrid orbits.
\end{theorem}

\begin{figure}
\begin{center}
\includegraphics[scale=.75]{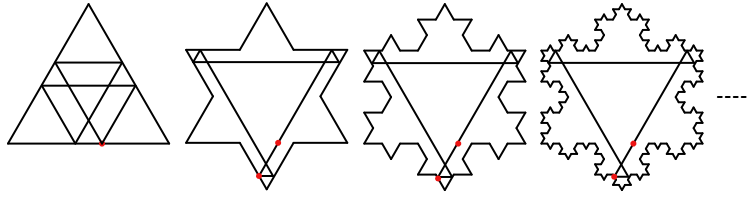}
\end{center}
\caption{An eventually constant sequence of compatible periodic hybrid orbits.  We see that the initial basepoint $\xoo=7/12$ lies on the middle third of the unit interval.  The basepoint $\xio{1}$ of the  compatible initial condition $(\xio{1},\frac{\pi}{3})$ has a ternary representation of type $\tern{lr}{\emptyset}$.}
\label{fig:CompatibleCantorOrbit7-12}
\end{figure}

\begin{example}[A constant sequence of compatible periodic hybrid orbits]
\label{exa:AConstantSequenceOfCompatiblePeriodicHybridOrbits}
Consider $\xoo=7/12$ in the base of the equilateral triangle. Such a value has a ternary representation of type $\tern{lr}{c}$.  Consider the initial condition $(\xoo,\frac{\pi}{3})$.  Then the sequence of compatible orbits $\compseqiang{1}{\frac{\pi}{3}}$ is a constant sequence. This follows from the fact that the ternary representation of $\xio{1}$ is $\overline{rl}$. Moreover, the representation of every basepoint of $\ofraciang{n}{\frac{\pi}{3}}$ is $\overline{lr}$.  In Figure \ref{fig:CompatibleCantorOrbit7-12}, we show the first three orbits in this (eventually) constant sequence of compatible periodic hybrid orbits.
\end{example}

As of now, the only examples of constant sequences of compatible periodic hybrid orbits are those for which the initial direction is $\frac{\pi}{3}$ and $\frac{\pi}{6}$ (and, equivalently, $\frac{\pi}{2}$).  When the initial angle of an orbit of a constant sequence of compatible periodic hybrid orbits is $\frac{\pi}{6}$ (or, equivalently, $\frac{\pi}{2}$), then the orbit will be degenerate.  For example, the orbit $\ofracixang{1}{\frac{3}{4}}{\frac{\pi}{2}}$ traverses a path that is a vertical line. This orbit has period $p=2$.  While $\compseqixang{1}{\frac{3}{4}}{\frac{\pi}{2}}$ is an important example of a constant sequence of compatible periodic hybrid orbits, it is arguably less interesting than the constant sequence of compatible periodic hybrid orbits $\compseqixang{1}{\frac{3}{4}}{\frac{\pi}{3}}$.

\vspace{2 mm}
\subsubsection{The corresponding prefractal translation surface $\sksi{n}$}
\label{subsubsec:TheCorrespondingPrefractalTranslationSurface}
\begin{figure}
\begin{center}
\includegraphics[scale=.3]{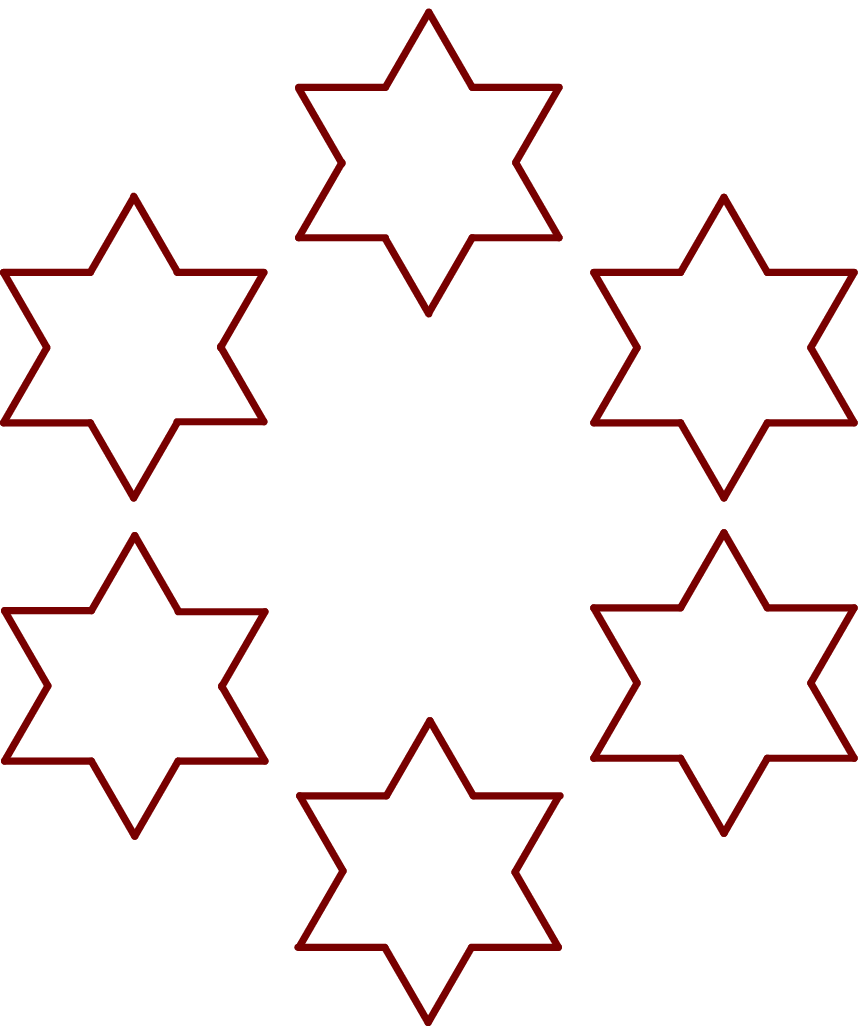}\quad\includegraphics[scale=.3]{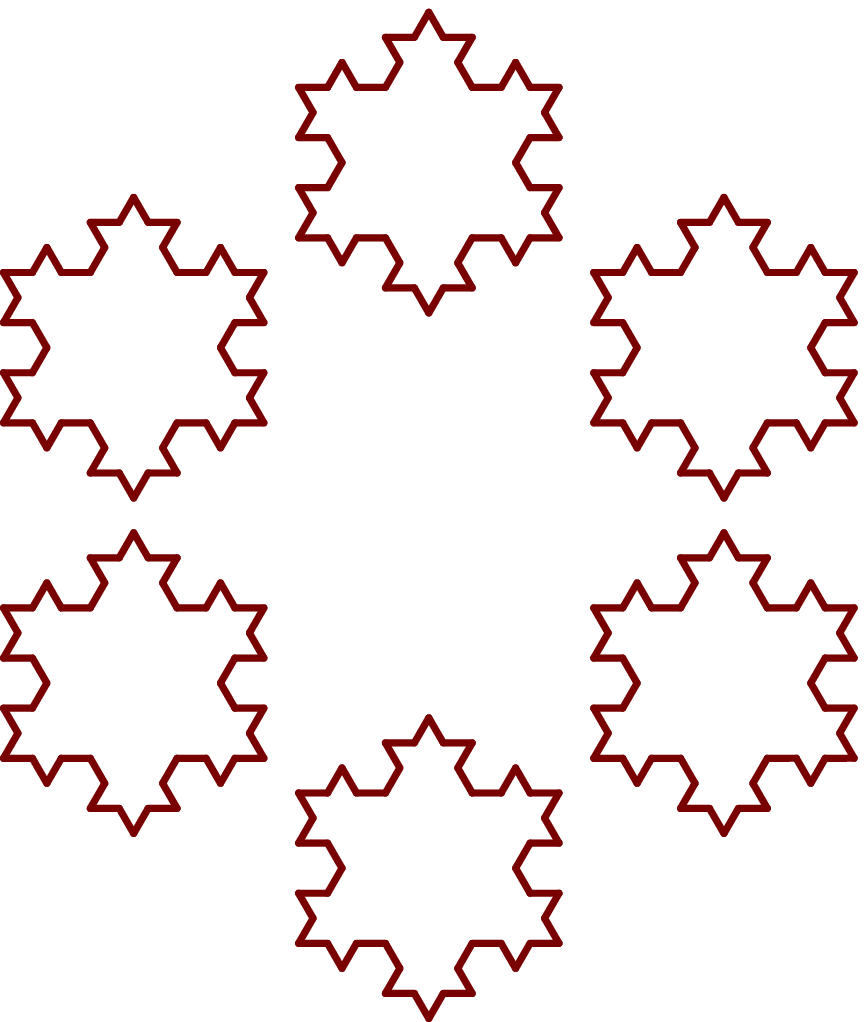}\quad\includegraphics[scale=.3]{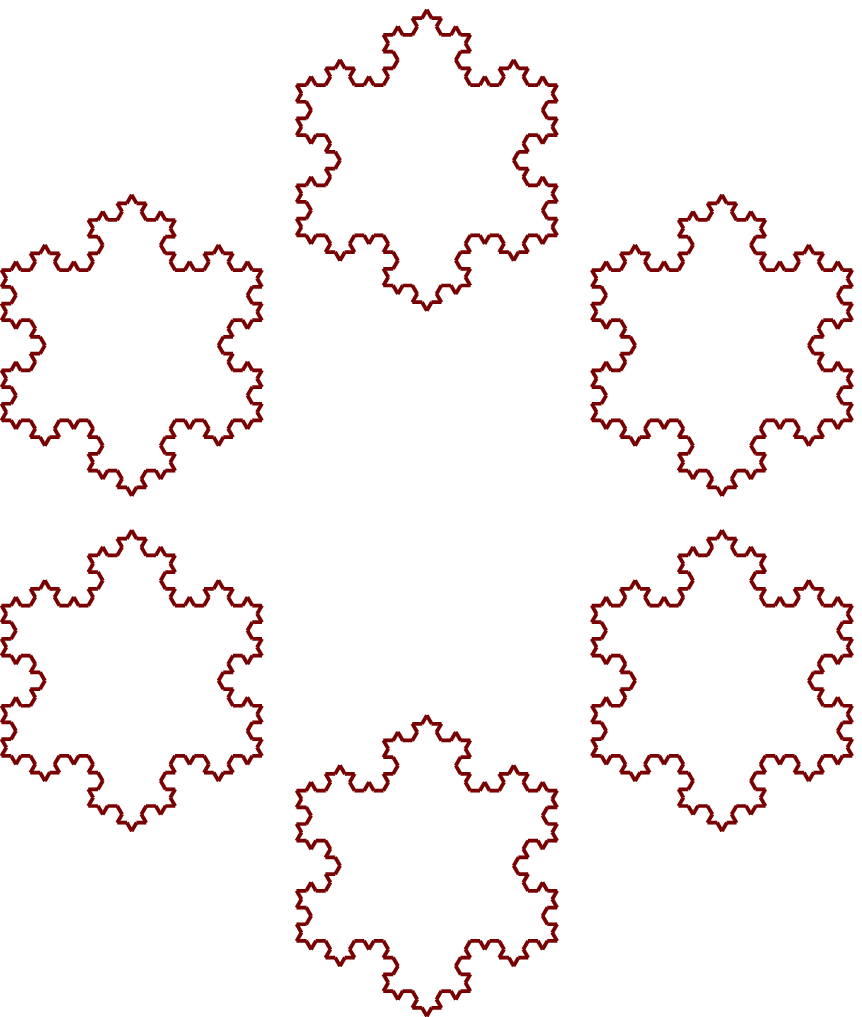}
\end{center}
\caption{The translation surfaces $\sksi{1}$, $\sksi{2}$ and $\sksi{3}$ associated with the Koch snowflake prefractal approximations $\ksi{1}$, $\ksi{2}$ and $\ksi{3}$, respectively.}
\label{fig:TheThreeCorrespondingFlatSurfaces}
\end{figure}

In \S\ref{subsec:translationStructuresandTranslationSurfaces} we saw how to construct a translation surface from a rational billiard table.  In the case of the equilateral triangle billiard table $\Omega(\Delta)=\omegaksi{0}$, there are $2\cdot \text{lcm}\{3,3,3\}=6$ copies of $\Omega(\Delta)$ used in the construction of the associated translation surface $\mathcal{S}(\Delta)$; see Example \ref{exa:EquilateralTriangleFlatSurface} and the associated Figure \ref{fig:sixEquiSurface}.  In the case of the prefractal billiard table $\omegaksi{n}$, only six copies of $\omegaksi{n}$ are needed in the construction of the associated translation surface $\mathcal{S}(\ksi{n})$, for every $n\geq 0$; see Figure \ref{fig:TheThreeCorrespondingFlatSurfaces}. (We refer to \cite{LapNie1,LapNie2,LapNie3} for further discussion of the topics in the present subsection.)

The vertices of $\omegaksi{n}$ correspond to conic singularities of the translation surface. However, only certain singularities are removable.  The vertices with angles measuring $\frac{\pi}{3}$ (measured from the interior), constitute removable singularities of the translation surface.  That is, the geodesic flow can be appropriately defined at these points.  The vertices with angles measuring $\frac{4\pi}{3}$ constitute nonremovable singularities.  Hence, it is possible to define reflection at certain vertices of the prefractal billiard $\omegaksi{n}$, but impossible to define at others.  Moreover, defining reflection at acute corners of $\omegaksi{n}$ in this way is independent of $n$.  That is, for a given vertex $v$ of $\omegaksi{n}$ with an acute angle $\frac{\pi}{3}$, the general rule for reflection in $v$ states that an incoming trajectory reflect through the angle bisector of $v$.  A billiard ball entering $v$ along the same path in $\omegaksi{n+1}$ as in $\omegaksi{n}$ will then reflect in $v$ in $\omegaksi{n+1}$ in exactly the same way as it did when considering $v$ as a vertex of $\omegaksi{n}$.

Such insight is clearly helpful in further understanding the behavior of a billiard ball on the Koch snowflake fractal billiard $\omegaks$, but we must be careful not to extrapolate more than is possible from this observation.  Knowing that we can determine an orbit of a prefractal billiard $\omegaksi{n}$ by unfolding the orbit of $\omegaksi{0}$ in $\omegaksi{n}$, we are inclined to allow orbits of $\omegaksi{0}$ that make collisions with corners.  However, a priori, we cannot conclude that such orbits do not unfold to form saddle connections in $\omegaksi{n}$ connecting two nonremovable singularities.  In the event an orbit $\ofraci{m}$ of $\omegaksi{m}$ intersects the boundary $\ksi{m}$ solely in acute corners, then such an orbit is an element of a sequence of compatible orbits $\compseqi{N}$ with $\ofraci{j}=\ofraci{m}$, for every $j\geq m$.

\subsection{The $T$-fractal prefractal billiard}
\label{subsec:TheTFractalPrefractalBilliard}

We refer to \S\ref{subsec:TheTFractal} for a discussion of the $T$-fractal $\mathscr{T}$ and of its prefractal approximations $\tfraci{n}$, for $n=0,1,2,...$; see, in particular, Figure \ref{fig:T-fractal}. Recall that the base of $\mathscr{T}_0$ has a length of two units. The prefractal billiard $\omegati{0}$ can be tiled by the unit square $Q$; see Figure \ref{fig:T0TiledByUnitSquare}. In general, for every $n\geq 0$, $\omegati{n}$ can be tiled by the square $\frac{1}{2^n}Q$.  As such, and since $Q$ obviously tiles the plane, we can apply Theorems \ref{thm:topologicalDichotomyForFn} and \ref{thm:generalTopologicalDichotomyForSequencesOfCompatibleOrbits}.

Much like the case of the prefractal Koch snowflake billiard $\omegaksi{n}$, we are interested in forming sequences of compatible orbits of prefractal billiards exhibiting particular properties.  The results in this subsection appear here for the first time and will be further discussed in \cite{LapNie6}. It is true that if a periodic orbit has an initial condition $(\xio{0},\theta_0^0)$, then there may exist a compatible orbit $\ofraci{N}$ that forms a saddle connection if $\xio{0}$ has a finite binary expansion.  This is not to suggest that $\ofraci{N}$ \textit{must} form a saddle connection. However, if every basepoint  $\xii{0}{k_0}$ of a periodic orbit $\mathscr{O}_Q(\xio{0},\theta_0^0)$ of the unit square has an infinite binary expansion (with no equivalent finite binary expansion), then viewing $\mathscr{O}_0(\xio{0},\theta_0^0)$ in $\omegati{0}$ as the reflected-unfolding of $\mathscr{O}_Q(\xio{0},\theta_0^0)$, the corresponding sequence of compatible orbits $\compseqi{0}$ will be a sequence of compatible periodic orbits.  We state this formally in the following theorem.

\begin{theorem}
Let $(\xio{0},\theta_0^0)$ be an initial condition of an orbit $\mathscr{O}_Q(\xio{0},\theta^0_0)$ of $\Omega(Q)$. Suppose every element of the footprint $\mathcal{F}_Q(\xio{0},\theta_0^0)$ has an infinite binary expansion \emph{(}and no equivalent finite binary expansion\emph{)} and $(\xio{0},\theta_0^0)$ is then the initial condition of an orbit of $\omegati{0}$ that constitutes the reflected-unfolding of $\mathscr{O}_Q(\xio{0},\theta_0^0)$ in $\omegati{0}$. Then the sequence of compatible orbits $\compseq$ \emph{(}where $(\xio{0},\theta^0) = (\xio{0},\theta_0^0)$\emph{)} of the prefractal billiards $\omegati{n}$ is a sequence of compatible periodic orbits.
\end{theorem}

\begin{figure}
\begin{center}
\includegraphics[scale=.4]{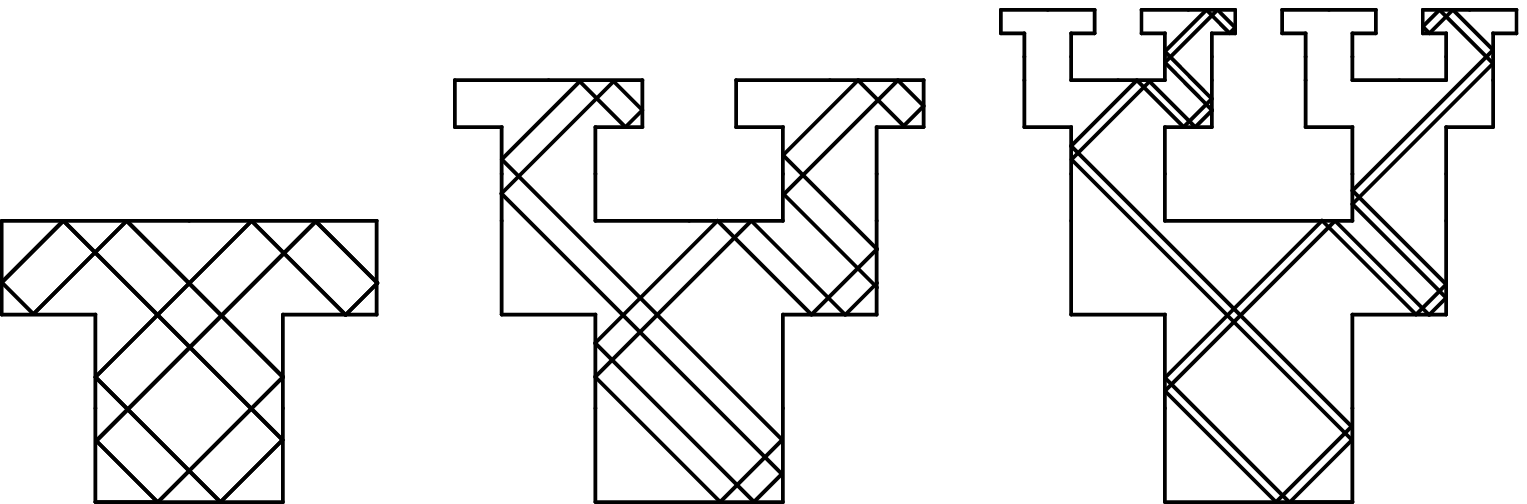}
\end{center}
\caption{A sequence of compatible periodic orbits of $\omegati{0}$, $\omegati{1}$ and $\omegati{2}$, respectively.}
\label{fig:T-FractalSequenceOfCompatibleOrbits}
\end{figure}

\begin{example}
\label{exa:T-fractalSequenceOfCompatibleOrbits}
Let $\xio{0} = \frac{4}{3}$ and $\theta_0^0 = \frac{\pi}{4}$.  Then, $\compseqang{\frac{\pi}{4}}$ is a nonconstant sequence of compatible periodic orbits; see Figure \ref{fig:T-FractalSequenceOfCompatibleOrbits}.
\end{example}

The following two theorems are ultimately concerning the prefractal billiard $\omegati{n}$.  Determining which intercepts and slopes yield line segments in the plane that avoid lattice points of the form $(\frac{a}{2^c},\frac{b}{2^d})$ is equivalent to specifying an initial condition of an orbit of a square billiard table that avoids corners of the billiard table.  Then, using the fact that an appropriately scaled square billiard table tiles $\omegati{n}$, we can reflect-unfold such an orbit in $\omegati{n}$ in order to determine an orbit of $\omegati{n}$.

\begin{theorem}
\label{thm:WhatTheSlopeCannotBeInTFrac}
Let $\xio{0} = \frac{t}{3^k}$ with $k,t\in \N$, $t$ and $3$ relatively prime, $k\neq 0$ and $0<t<3^k$.  Further, let $m\in\R$.  If for every $p,q,r,s\in \Z$, $r,s\geq 0$, we have that
\begin{align}
m &\neq \frac{q2^{r-s}3^k}{p3^k-t2^r},
\label{eqn:WhatmCannotBe}
\end{align}
\noindent then the line $y=m(x-\xio{0})$ does not contain any point of the form $(\frac{a}{2^c},\frac{b}{2^d})$, $a,b,c,d\in \Z$, with $c,d\geq 0$.
\end{theorem}

Note that the condition (\ref{eqn:WhatmCannotBe}) above is automatically satisfied if the slope $m$ is irrational.


\begin{theorem}
\label{thm:WhatTheSlopeCanBeInTFrac}
Let $\xio{0} = \frac{t}{3^k}$, with $k,t\in\N$, $t$ and $3$ relatively prime, $k\neq 0$ and $0<t<3^k$.  If
\begin{align}
\notag m&=\frac{2^\gamma}{(2\alpha +1)^\beta},
\end{align}
\noindent with $\alpha,\beta,\gamma\in \N$, $\alpha,\beta,\gamma\geq 0$, then, for every $p,q,r,s\in \Z$ with $r,s\geq 0$, the point $(\frac{p}{2^r},\frac{q}{2^s})$ does not lie on the line $y=m(x-\xio{0})$.
\end{theorem}


Finally, Theorems \ref{thm:WhatTheSlopeCannotBeInTFrac} and \ref{thm:WhatTheSlopeCanBeInTFrac} combined with the fact that an initial condition of an orbit of $\omegati{N}$, $N\geq 0$, determines a sequence of compatible orbits $\compseqi{N}$, allows us to determine a countably infinite family of sequences of compatible periodic orbits.

\vspace{2 mm}
\subsubsection{The corresponding prefractal translation surface $\mathcal{S}(\mathscr{T}_n)$}

\begin{figure}
\begin{center}
\includegraphics[scale=.6]{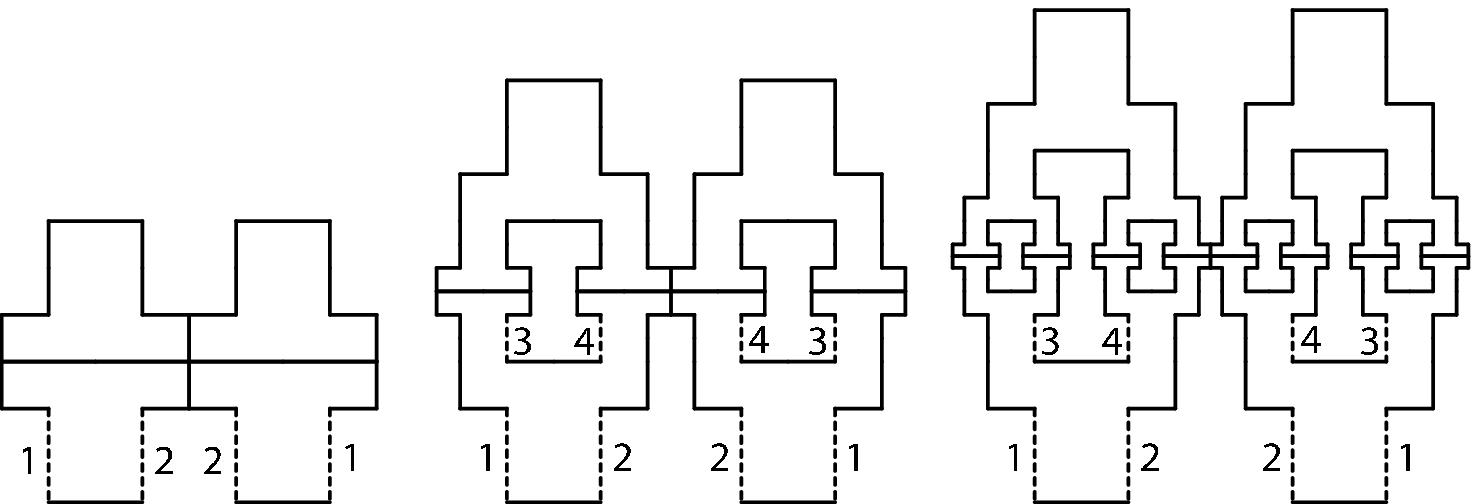}
\end{center}
\caption{The translation surfaces $\stfraci{0}$, $\stfraci{1}$ and $\stfraci{2}$ associated with the $T$-fractal prefractal approximations $\tfraci{0}$, $\tfraci{1}$ and $\tfraci{2}$, \mbox{respectively.}}
\label{fig:T-fractalFlatSurfaces}
\end{figure}

For every $n\geq 0$, the interior angles of $\tfraci{n}$ are $\frac{\pi}{2}$ and $\frac{3\pi}{2}$.  To form the associated translation surface $\stfraci{n}$, we appropriately identify four copies of $\omegati{n}$; see Figure \ref{fig:T-fractalFlatSurfaces} for a depiction of the first three translation surfaces.

Then, every point of $\stfraci{n}$ associated with a vertex of $\omegati{n}$ measuring $\frac{\pi}{2}$ constitutes a removable singularity of $\stfraci{n}$.  Similarly, every point of $\stfraci{n}$ associated with a vertex of $\omegati{n}$ of interior angle measuring $\frac{3\pi}{2}$ constitutes a nonremovable singularity of $\stfraci{n}$.  Therefore, not every vertex of $\omegati{n}$ will present a problem for the billiard flow.

Consider an orbit of $\Omega(Q)$, where the orbit has basepoints corresponding to vertices of $Q$, the unit square.  Since such vertices correspond to removable singularities in the corresponding translation surface (this being the flat torus, see \S\ref{subsec:translationStructuresandTranslationSurfaces}), we see that the same orbit reflected-unfolded in the billiard $\omegati{n}$ (if one first scales the billiard $\Omega(Q)$ and the orbit contained therein by $\frac{1}{2^n}$, see \S\ref{subsec:UnfoldingABilliardOrbit}) can potentially intersect vertices of $\tfraci{n}$ that are associated with nonremovable singularities in the corresponding translation surface.


\subsection{A prefractal self-similar Sierpinski carpet billiard}
\label{subsec:AprefractalSelfSimilarSierpinskiCarpetBilliard}
Let $S_{a}$ be a self-similar Sierpinski carpet, as defined in Definition \ref{def:ASelfSimilarSierpinskiCarpet}, and let us denote its natural prefractal approximations by $\sai{i}$ for $i=0,1,2,...$ (as in \S\ref{subsec:aSierpinskiCarpet}).  The corresponding billiard is then denoted by $\omegasa{a}$.  In this subsection, we examine the behavior of the billiard flow on the rational polygonal billiard given by the prefractal approximations $\omegasai{a}{i}$.\footnote{We note that the results in this subsection appear here for the first time and will be further discussed in \cite{CheNie}.}  In the event a billiard ball collides with a corner of a peripheral square, we must terminate the flow and such a trajectory is then called \textit{singular}.  In addition to being singular, such a trajectory will form a saddle connection (see the beginning of \S\ref{sec:RationalBilliards} for a discussion of closed billiard orbits that form saddle connections). As we have discussed, an examination of the corresponding translation surface may prove useful in determining whether or not a billiard ball can reflect in a vertex.

\begin{definition}[Obstacle of $\Omega(D)$]
\label{def:ObstacleToBilliardFlow}
Let $\Omega(D)$ be a polygonal billiard. Then $\Omega(D)$ can be modified by placing in its interior a piecewise smooth segment that inhibits the billiard flow and causes a billiard ball to reflect.  Such a segment is called an \textit{obstacle} of $\Omega(D)$.
\end{definition}

Clearly, each prefractal billiard $\omegasai{a}{i}$ can be interpreted as a \textit{square billiard with obstacles}.


\begin{notation}
\label{nota:AlphaForTheta}
Due to the fact that Theorem \ref{thm:AsetBset} refers to the slope of a nontrivial line segment and we make heavy use of this theorem, we will denote the initial condition $(\xio{n},\theta_n^0)$ of an orbit of $\omegasi{n}$ by $(\xio{n},\alphaio{n})$, where $\alphaio{n} = \tan(\theta_n^0)$.
\end{notation}

\begin{definition}[An orbit of the cell $\celli{k}{a^k}$ of $\omegasi{k}$]
 Consider the boundary of a cell $\celli{k}{a^k}$ of $\omegasi{k}$ as a barrier.\footnote{Here, $C_{k,a^k}$ is a cell of the $k$th prefractal approximation $\sai{k}$, as given in Definition \ref{def:ACellOfSai} with all numbers $a_j$ equal to $a$.}  Then an orbit with an initial condition contained in the cell is called an \textit{orbit of the cell $\celli{k}{a^k}$ of $\omegasi{k}$}.
\end{definition}

\begin{remark}
So as to be clear, the boundary of the cell does not form an obstacle to the billiard flow, as defined in Definition \ref{def:ObstacleToBilliardFlow}.  Rather, we are treating the cell $\celli{k}{a^k}$ as a billiard table in its own right, embedded in the larger prefractal approximation $\omegasi{k}$.
\end{remark}

Recall from \S\ref{subsec:aSierpinskiCarpet} that a self-similar Sierpinski carpet $S_a$ is the unique fixed point attractor of a suitably chosen iterated function system $\{\phi_{j}\}_{j=1}^{a^2-1}$ consisting of similarity contractions.  In light of this, an orbit of a cell $\celli{k}{a^k}$ of $\omegasi{k}$ is the image of an orbit $\osi{0}$ of the unit-square billiard $\Omega(S_0)$ under the action of a composition of contraction mappings $\phi_{m_k}\circ\cdots\circ\phi_{m_1}$, with $1\leq m_i\leq a^2-1$ and $1\leq i\leq k$, determined from the iterated function system $\{\phi_j\}_{j=1}^{a^2-1}$ of which $S_a$ is the unique fixed point attractor.

\begin{lemma}
\label{lem:SegmentBeginningAtMidpointIsNonTrivialIfAlphaFromB}
Consider a self-similar Sierpinski carpet $S_a$. Let $k\geq 0$ and $S_{a,k}$ be a prefractal approximation of $S_a$. If $\alpha\in B_a$,\footnote{Recall from Notation \ref{not:AaBaAbBb} that $B_a$ is the set of slopes given by Equation (\ref{eqn:Bset}).} then the line segment beginning at a midpoint of a cell $\celli{k}{a^k}$ of $S_{a,k}$ is a nontrivial line segment \emph{(}in the sense of Definition \ref{def:nontrivialLineSegmentOfSa}\emph{)}.  Moreover, such a segment avoids the boundary of the peripheral squares of $S_a$ with side-length $a^{-m}$, $m\geq k+1$.
\end{lemma}


The statement in Lemma \ref{lem:SegmentBeginningAtMidpointIsNonTrivialIfAlphaFromB} asserts that a segment beginning at a midpoint of a cell with slope $\alpha\in B_a$ will be a nontrivial line segment in $S_a$.  In addition to this, any line segment contained in $\mathbb{R}^2$ that contains a nontrivial line segment of $S_a$ must necessarily avoid the peripheral squares in a tiling of $\mathbb{R}^2$ by $S_a$.  Otherwise, there exists $k\geq 1$ such that scaling the line segment in $\mathbb{R}^2$ and the tiling of $\mathbb{R}^2$ by $a^{-k}$ results in a segment contained in the nontrivial line segment which intersects peripheral squares of $S_a$.  This is a contradiction of the fact that the segment beginning at $(2^{-1},0)$ with slope $\alpha\in B_a$ is a nontrivial line segment of $S_a$.  We then deduce the following result.

\begin{theorem}
\label{thm:ConstantSequenceOfCompatibleOrbitsInSa}
Consider a self-similar Sierpinski carpet $S_a$.  Let $k\geq 0$ and $S_{a,k}$ be a prefractal approximation of $S_a$.  Furthermore, let $\alpha\in B_a$ and $\xio{k}=(p(2a^k)^{-1},0)$ with $p\leq a^k$ a positive, odd integer. If $\osi{k}$ is an orbit of $\omegasai{a}{k}$, then the initial condition $(\xio{k},\alphaio{k})$ determines a sequence of compatible periodic orbits $\compseqsi{k}$ of the prefractal approximations $\omegasai{a}{n}$.
\end{theorem}

As one may suspect, there exists $N\geq k\geq 0$ such that a sequence of compatible orbits $\compseqsi{N}$ is a constant sequence of compatible orbits.  Moreover, $\xio{n}=\xio{N}$, for every $n\geq N$.  This is not any different from the case of a constant sequence of compatible orbits of prefractal billiards $\omegaksi{n}$, as discussed in Theorem \ref{thm:SufficientConditionForCantorOrbit} and Example \ref{exa:AConstantSequenceOfCompatiblePeriodicHybridOrbits}.  However, in the context of a self-similar Sierpinski carpet billiard table, every sequence of compatible orbits we will examine will be a sequence for which there exists $N\geq 0$ such that $\compseqsi{N}$ is a constant sequence of compatible orbits.

\vspace{2 mm}
\subsubsection{The corresponding prefractal translation surface $\ssai{i}$}

\begin{figure}
\begin{center}
\includegraphics[scale=.2]{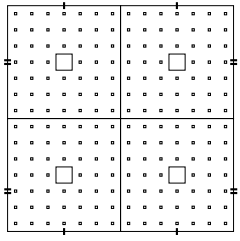}
\end{center}
\caption{Interpreting the translation surface $\ssai{n}$ as a flat torus with obstacles.}
\label{fig:squareTorusWithObstacles}
\end{figure}

In much the same way the billiard $\omegasai{a}{i}$ can be interpreted as a square billiard with obstacles, the corresponding translation surface can be interpreted as a ``torus with obstacles''; see Figure \ref{fig:squareTorusWithObstacles}.

In light of the fact that $\ssai{n}$ can be interpreted as a torus with obstacles and the presence of a  dynamical equivalence between the billiard flow and the geodesic flow on the corresponding translation surface (see \S\ref{subsec:UnfoldingABilliardOrbit}), we see that reflection in the vertices with angles measuring $\frac{\pi}{2}$ (relative to the interior) can be defined.  More specifically, the geodesic flow can be defined at points corresponding to vertices with angles measuring $\frac{\pi}{2}$, because these points constitute removable singularities of the geodesic flow.

This fact is crucial in determining orbits of $\omegasa{a}$ for which the slope $\alpha$ is an element of $A_a$ and not $B_a$ (see Notation \ref{not:AaBaAbBb}), and the orbit avoids all peripheral squares of $\omegasa{a}$.  While one may say that this contradicts part of Theorem \ref{thm:AsetBset} (and he/she would be right), in \cite{CheNie} a more precise formulation of Theorem \ref{thm:AsetBset} is given that clarifies which slopes are permissible and which ones are not.  That is, if $\alphaio{n}\in A_a$, it may be possible for an orbit $\osi{n}$ to begin at the origin and avoid the peripheral squares of each billiard $\omegasai{a}{m}$, for every $m\geq n$.  We do not give here an explicit reformulation of Theorem \ref{thm:AsetBset}, but Example \ref{exa:OrbitExampleContradictingTysonsTheoremSlope2-3OriginS7} in \S\ref{subsec:aSelf-SimilarSierpinskiCarpetBilliard} exhibits a situation showing that the latter half of Theorem \ref{thm:AsetBset} is not stated precisely enough.

\section{Fractal billiards}
\label{sec:FractalBilliards}

The theme that will tie together all of the examples in \S\ref{sec:PrefractalRationalBilliards} is that suitable limits of sequences of compatible orbits may constitute billiard orbits of each respective fractal billiard table.  We have shown that in the case of $\omegaks$, $\omegat$ and $\omegasa{a}$, we can determine a sequence of compatible periodic orbits.  We will see that in each case of a fractal billiard, under certain conditions, a sequence of compatible periodic orbits (or a proper subset of points from each footprint $\fprintfraci{n}$) will converge to a set which can be thought of as a true orbit of a fractal billiard table (or such a sequence will yield a subsequence of basepoints converging to what we are calling an \textit{elusive point} in \cite{LapNie2,LapNie3}).

\subsection{A general framework for $\omegaks$, $\omegat$ and $\Omega(S_a)$}
\label{subsec:aGeneralStructureForFractalBilliards}
We restrict our attention to the family of fractal billiard tables $\omegafrac$ where $F$ is a fractal approximated by a suitable sequence of rational polygons $\{F_n\}_{n=0}^\infty$, with each $F_n$ tiled by $D_n=c_nD_0$ for suitably chosen $c_n\in(0,1]$ and $D_0$ a polygon that tiles the plane.  Specifically, we are interested in developing a general framework for dealing with a fractal billiard table $\Omega(F)$ which   is similar to that of $\omegaks$, $\omegat$ and $\Omega(S_a)$.

Before we begin our discussion of the fractal billiard tables $\omegaks$, $\omegat$ and $\omegasa{a}$, we define certain terms.  The following definitions were initially motivated by the work in \cite{LapNie3}, but later generalized for this paper in order to account for a larger class of fractal billiard tables. (From now on, we assume that $\Omega(F)$ is a fractal billiard table with prefractal billiard approximations $\{\omegafraci{n}\}_{n=0}^\infty$ as described just above.)

\begin{definition}[A corner]
\label{def:ACorner}
Let $z\in F$. If there exists $n\geq 0$ such that $z\in F_n$ and $z$ is a vertex of $F_n$, then $z$ is called a \textit{corner} of $F$.
\end{definition}

\begin{definition}[A Cantor point]
\label{def:ACantorPoint}
Let $z\in F$ be such that $z$ is not a corner of $F$. If there exists $N\geq 0$ such that for every $n\geq N$, $z\in F_n$ and every connected neighborhood of $z$ contained in $F_n$ becomes totally disconnected when intersected with $F$, then $z$ is called a \textit{Cantor point} of $F$.
\end{definition}

In the Koch snowflake $\ks$, every Cantor point is a smooth point of infinitely many prefractals $\ksi{n}$ approximating $\ks$.  That is, if $z$ is a Cantor point in $\ks$, then there exists $N\geq0$ such that for every $n\geq N$, there exists a well-defined tangent at $z\in \ksi{n}$.\footnote{Here and Definition \ref{def:SmoothFractalPoint} below, $z$ is viewed as a point of the smooth subarc of $F_n$ to which it belongs.}  We deduce from this that the law of reflection holds at $z\in \ksi{n}$, for every $n\geq N$.  Moreover, since the billiard ball reflects at $z\in \ksi{n}$ at the same angle for every $n\geq N$, we deduce that the tangent at $z$ is the same for each $\ksi{n}$, $n\geq N$. This observation then prompts us to generalize the definition of a Cantor point in order to account (for example) for points of the $T$-fractal which are not Cantor points, but are points for which a well-defined tangent can be found in infinitely many prefractal approximations.

\begin{definition}[Smooth fractal point]
\label{def:SmoothFractalPoint}
Let $z\in F$ and $N\geq 0$ be such that $z\in F_n$ for every $n\geq N$.  If there exists a well-defined tangent at $z\in F_n$ for every $n\geq N$, then $z$ is called a \textit{smooth fractal point}.
\end{definition}

To be clear, a Cantor point of $F$ is an example of a smooth fractal point of $F$. The special nature of a Cantor point warrants a formal definition. In the $T$-fractal billiard, there are certainly corners and elusive points.  There are also smooth connected segments contained in the boundary of $\omegat$.  Points contained in such segments that do not correspond to corners are then called smooth fractal points.

\begin{definition}[An elusive point]
\label{def:AnElusivePoint}
Let $z\in F$. If $z\notin \bigcup_{n=0}^\infty F_n$, then $z$ is called an \textit{elusive point} of $F$.
\end{definition}

Consider a piecewise linear path in $\omegafrac$, such that every linear segment of the path is joined at the endpoint of another segment with the coincidental endpoints intersecting the boundary $F$ at a smooth fractal point of $F$ (in the sense of Definition \ref{def:SmoothFractalPoint}).  In the following definition, we define a particular type of piecewise linear curve in a fractal billiard $\omegafrac$.

\begin{definition}[A nontrivial path]
Suppose that there exists a piecewise linear curve in $\omegafrac$ as described immediately above.  If at each point $z$ for which the piecewise linear path intersects the boundary $F$, the angle formed by the first segment is equal to the angle formed by the second segment, relative to the side of $F_n$ on which $z$ lies,\footnote{Recall from Definition \ref{def:SmoothFractalPoint} that a smooth fractal point $z$ of $F$ is necessarily a point of infinitely many prefractal approximations $F_m$.  Hence, there is a least nonnegative integer $n$ such that $z\in F_n$.} then the piecewise linear path is called a \textit{nontrivial path} of $\omegafrac$.
\end{definition}

\begin{remark}
In \cite{LapNie3}, a nontrivial path was called a \textit{nontrivial polygonal path}.  The change in name is purely based on aesthetics.
\end{remark}

\begin{definition}[A Cantor orbit]
\label{def:ACantorOrbit}
Suppose $\ofraci{N}$ is an orbit of $\omegafraci{N}$, for some $N\geq 0$, such that every point of the footprint $\fprintfraci{N}$ corresponds to a smooth fractal point of $F$.  This then readily implies that $\ofraci{n}$ is the same as $\ofraci{N}$ for every $n\geq N$.\footnote{In other words, $\compseqi{N}$ is a constant sequence of compatible orbits, where a \textit{sequence of compatible orbits} was defined in Definition \ref{def:SequenceOfCompatibleOrbits}.} Then $\ofraci{n}$ is called a \textit{Cantor orbit} of $\omegafrac$ and is denoted by $\mathscr{O}(x^0,\theta^0)$.
\end{definition}

If $\omegafrac$ is a fractal billiard table, then it may or may not be possible to construct Cantor orbits or nontrivial paths of $\Omega(F)$.  We will next discuss three examples of fractal billiard tables with different dynamical properties that lend themselves well (or not) to determining well-defined billiard orbits.

\begin{remark}
We note that applying Definitions \ref{def:ACorner}, \ref{def:ACantorPoint} and \ref{def:AnElusivePoint} to $\omegaks$ and the sequence of rational polygon prefractal approximations $\omegaksi{n}$ which we have discussed in \S\ref{subsec:TheKochCurveAndKochSnowflake} yields exactly the sets of points we are considering as corners, Cantor points and elusive points of $\omegaks$, respectively.  Moreover, applying Definitions \ref{def:ACorner}, \ref{def:SmoothFractalPoint} and \ref{def:AnElusivePoint} to $\omegat$ and the prefractal approximations $\omegati{n}$ which we discussed in \S\ref{subsec:TheTFractal} yields exactly the sets of points that we are considering as corners, smooth fractal points and elusive points of $\omegat$.  Finally, applying Definitions \ref{def:ACorner} and \ref{def:SmoothFractalPoint} to $\Omega(S_a)$ and the prefractal approximations $\omegasai{a}{n}$ which we discussed in \S\ref{subsec:aSierpinskiCarpet} yields exactly the set of points we are considering as corners and smooth fractal points of $\Omega(S_a)$.
\end{remark}

\subsection{The Koch snowflake fractal billiard}
\label{subsec:TheKochSnowflakeFractalBilliard}
As we have noted before at the end of \S\ref{subsec:TheKochCurveAndKochSnowflake}, for each $n\geq 0$, $\ksi{n}\cap \ks$ can be realized as the union of $3\cdot 4^n$ self-similar ternary Cantor sets, each spanning a distance of $\frac{1}{3^n}$.  Within each Cantor set, we find Cantor points and corners of the Koch snowflake.

We begin our discussion of orbits of $\omegaks$ by examining the limiting behavior of a particular sequence of compatible orbits with the initial condition $(\xio{N},\frac{\pi}{3})$, where $\xio{N}$ is a Cantor point of $\ks$ (i.e., $\xio{N}$ is a point of $\ksi{N}$ with a well-defined tangent
 in $\ksi{n}$ for every $n\geq N$).  For the sake of simplicity, we let $N=0$ and $\xio{N} = \frac{1}{4}$ be on the base of the equilateral triangle $\ksi{0}$ (recall that we are assuming that the left corner of $\ksi{0}$ is at the origin and the length of each side is one unit).  Then, $\ofraciang{0}{\frac{\pi}{3}}$ is an orbit that remains fixed as one constructs $\omegaksi{1}$ from $\omegaksi{0}$.  More correctly, $\{\ofraciang{n}{\frac{\pi}{3}}\}_{n=0}^\infty$ is a sequence of compatible orbits with $\fprintfraciang{n}{\frac{\pi}{3}}= \fprintfraciang{0}{\frac{\pi}{3}}$ for every $n\geq 0$ (that is, with the same footprint in each prefractal approximation).

In general, if $(\xio{N},\frac{\pi}{3})$ is an initial condition of an orbit of $\omegaksi{N}$ and $\xio{N}$ is a Cantor point, then the sequence of compatible orbits is such that for every $n\geq N$, the footprints $\fprintfraciang{n}{\frac{\pi}{3}}$ and $\fprintfraciang{N}{\frac{\pi}{3}}$ are the same.

\begin{theorem}
If $x\in \ks$ is a Cantor point, then there exists a well-defined orbit of $\omegaks$ with an initial condition $(x,\varpi(\frac{\pi}{3}))$, where the angle $\varpi(\frac{\pi}{3})$ is determined with respect to the side on which $x$ lies in a prefractal approximation $\omegaksi{n}$.
\end{theorem}





There are many more properties of $\compseqiang{N}{\frac{\pi}{3}}$ which we could discuss here.  These properties largely rely on the nature of the ternary representation of $\xio{N}$, and are elaborated upon in \cite{LapNie2,LapNie3}.  We now proceed to illustrate how we can connect two elusive points of $\omegaks$.  Such a result has already been presented in greater detail in \cite{LapNie3}, so we will be brief.   In \S\ref{subsec:TheTFractalBilliard}, we will show that an identical construction holds for the billiard table $\omegat$.

\begin{figure}
\begin{center}
\includegraphics[scale=.3]{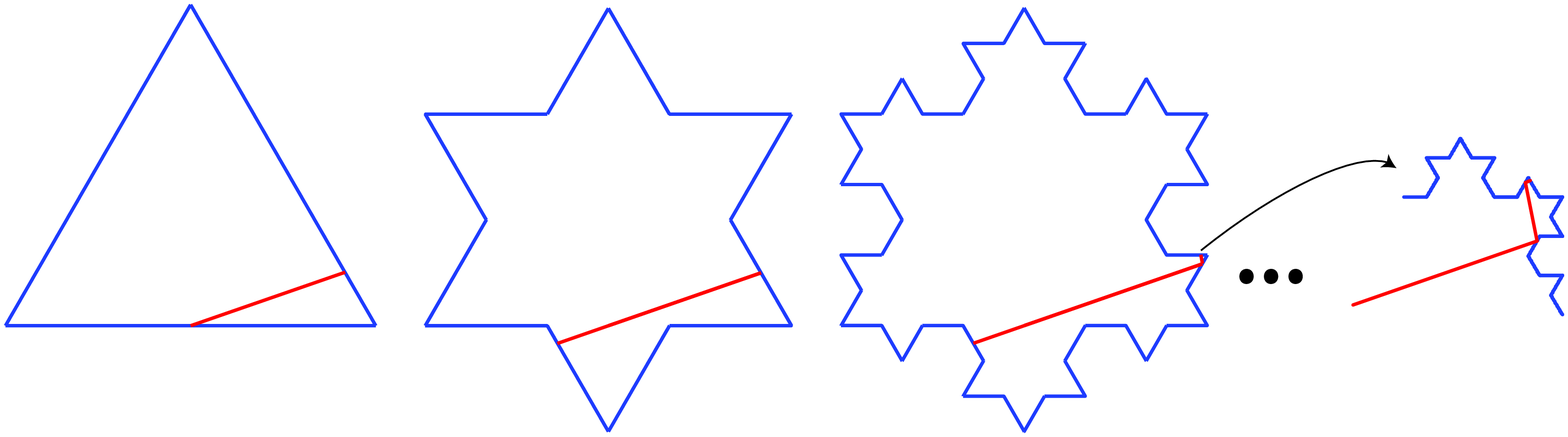}
\end{center}
\caption{A nontrivial path of the Koch snowflake fractal billiard table $\omegaks$ beginning at $x=\frac{1}{2}$.}
\label{fig:nontrivialPath}
\end{figure}

\begin{figure}
\begin{center}
\includegraphics[scale=.5]{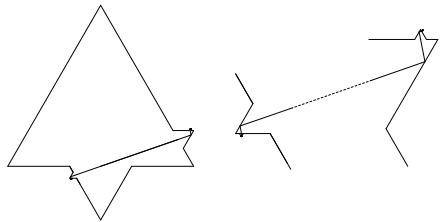}
\end{center}
\caption{Two nontrivial paths connecting two elusive points of $\omegaks$. (As is explained in the text, these two paths can be concatenated to obtain a single nontrivial path connecting the two elusive points.)  In the first figure, we only show the relevant portions of the Koch snowflake.  In the second figure, we magnify the regions containing the nontrivial paths so as to highlight the fact that such paths are converging to elusive points.  Actually, there is an obvious geometric similarity one can take advantage of in order to produce more segments of the nontrivial path.}
\label{fig:twoNontrivialPaths}
\end{figure}

Recall from Example \ref{exa:ASequenceOfCompatiblePeriodicHybridOrbits} that we were able to construct a sequence of compatible periodic hybrid orbits.  From such a sequence we can derive a sequence of basepoints that is converging to an elusive point of $\omegaks$.  The latter sequence of basepoints constitutes the vertices of a nontrivial path; see Figure \ref{fig:nontrivialPath}.  One may consider a direction $\gamma_0^0$ that is the reflection of $\theta_0^0$ through the normal at $\xio{0}$.  Then, the resulting sequence of compatible periodic hybrid orbits $\compseqiang{N}{\gamma_0^0}$ yields a sequence of basepoints converging to another elusive point.  Again, such a sequence of basepoints constitutes the vertices of a nontrivial path of $\omegaks$; see Figure \ref{fig:twoNontrivialPaths}.  Together, these two nontrivial paths constitute a single nontrivial path connecting two elusive points of $\omegaks$.

In conjunction with Theorem \ref{thm:bodd}, we can determine countably infinitely many initial conditions $(\xio{n},\theta_n^0)$, each of which determines a sequence of compatible periodic hybrid orbits yielding a sequence of basepoints converging to an elusive point of $\omegaks$.

\subsection{The $T$-fractal billiard}
\label{subsec:TheTFractalBilliard}
\begin{figure}
\begin{center}
\includegraphics[scale=.5]{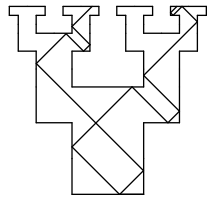}
\end{center}
\caption{Two nontrivial paths connecting two elusive points of $\omegat$.}
\label{fig:T-fractalNontrivialPaths}
\end{figure}

The results in this subsection appear here for the first time and will be further discussed in \cite{LapNie6}. We begin our discussion of the billiard $\omegat$ by recalling (and referring the reader back to) Example \ref{exa:T-fractalSequenceOfCompatibleOrbits} from \S\ref{subsec:TheTFractalPrefractalBilliard}.  The sequence of compatible periodic orbits provided by Example \ref{exa:T-fractalSequenceOfCompatibleOrbits} gives rise to a nontrivial path that connects $\frac{4}{3}$ with an elusive point of $\omegat$.  Furthermore, considering the sequence of compatible periodic orbits $\compseqixang{N}{\frac{4}{3}}{\frac{3\pi}{4}}$, we determine another nontrivial path that connects $\frac{4}{3}$ with another elusive point of $\omegat$; see Figure \ref{fig:T-fractalNontrivialPaths}.  This behavior is analogous to the one which we observed for the Koch snowflake billiard in \S\ref{subsec:TheKochSnowflakeFractalBilliard}.

As was the case with $\omegaks$, we can analogously build upon Theorems \ref{thm:WhatTheSlopeCannotBeInTFrac} and \ref{thm:WhatTheSlopeCanBeInTFrac} in order to determine a sequence of basepoints converging to an elusive point. That is, Theorems \ref{thm:WhatTheSlopeCannotBeInTFrac} and \ref{thm:WhatTheSlopeCanBeInTFrac} guide our search for a sequence of compatible periodic orbits which yields a sequence of basepoints converging to an elusive point of $\omegat$.

\begin{theorem}
Let $\compseqi{N}$ be a sequence of compatible orbits.  Then, there are countably infinitely many directions and countably infinitely many points from which to choose so that $\compseqi{N}$ is a sequence of compatible periodic orbits yielding a sequence of basepoints $\{x_{n}^{k_n}\}_{n=N}^\infty$ that converges to an elusive point of $\omegat$. The collection of basepoints $\{x_{n}^{k_n}\}_{n=N}^\infty$ constitute the vertices of a nontrivial path of $\omegat$. Moreover, once such a nontrivial path is constructed, letting ${x^\prime}_N^0 = x_N^0$,  an additional nontrivial path can be determined from a sequence of compatible periodic orbits $\compseqixang{N}{{x^\prime}_n^0}{\pi-\theta^0}$ in exactly the same fashion.
\end{theorem}


\subsection{A self-similar Sierpinski carpet billiard}
\label{subsec:aSelf-SimilarSierpinskiCarpetBilliard}

In \cite{Du-CaTy}, nontrivial line segments of Sierpinski carpets are constructed.  Building on the main results of \cite{Du-CaTy}, the second author and Joe P. Chen have been able to construct a family of Cantor periodic orbits of a self-similar Sierpinski carpet, in the sense of \cite{LapNie2,LapNie3} recalled in Definition \ref{def:ACantorOrbit}.\footnote{The results described in this subsection appear here for the first time and will be further discussed in \cite{CheNie}.}  Such orbits constitute Cantor orbits of the self-similar Sierpinski carpet.  As of yet, we have not attempted to construct a nontrivial path of a Sierpinski carpet.

In light of Theorem \ref{thm:ConstantSequenceOfCompatibleOrbitsInSa}, we say that the trivial limit of a constant sequence of compatible periodic orbits constitutes a periodic orbit of a self-similar Sierpinski carpet billiard $\omegasa{a}$.  In the event an orbit has an initial direction $\alphaio{0}$, we may still be able to determine a constant sequence of compatible periodic orbits.  The trivial limit of such a sequence then constitutes a periodic orbit of $\omegasa{a}$. Using the fact that reflection can be defined in the vertices with interior angles measuring $\frac{\pi}{2}$, we can state the following result. (Recall from \S\ref{subsec:AprefractalSelfSimilarSierpinskiCarpetBilliard} that $S_{a,n}$ is the $n$th prefractal approximation of $S_a$.)

\begin{theorem}
\label{thm:AlphaInAFormsOrbitOfSa}
Recall from Notation \ref{nota:AlphaForTheta} that if $\theta$ is the initial direction of a billiard orbit, then $\alpha = \tan\theta$.  Let $x^0=(0,0)$, $\alpha\in \mathbb{Q}$ and let $\mathscr{O}(x^0, \alpha)$ be an orbit of $\Omega(S_0)$.  If $\mathscr{O}(x^0,\alpha)$, as an orbit of $\omegasai{a}{1}$, avoids the middle peripheral square, then the initial condition $(x^0,\alpha)$ will determine an orbit of $\omegas$.  Specifically, the path traversed by the orbit $\mathscr{O}(x^0,\alpha)$ of $\omegasai{a}{1}$ is exactly the path traversed by the orbit of $\omegasa{a}$ determined by $(x^0,\alpha)$.
\end{theorem}

\begin{figure}
\begin{center}
\includegraphics[scale=.4]{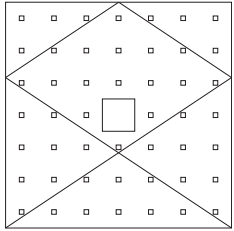}
\end{center}
\caption{An orbit with an initial condition beginning at $(0,0)$ and with an initial direction constituting a slope of $\alpha=2/3\in A_7$, where $A_7$ is defined as in Notation \ref{not:AaBaAbBb}.  While it would appear that this orbit intersects corners of peripheral squares, it in fact remains away from all peripheral squares.  The same is true for finer approximations.}
\label{fig:2-3orbitStartingFrom1-2}
\end{figure}

\begin{example}
\label{exa:OrbitExampleContradictingTysonsTheoremSlope2-3OriginS7}
Let $x^0=(0,0)$, $\alpha=2/3\in \slopesa{5}$. Consider an orbit of $\Omega(S_{7,2})$ with an initial condition $(x^0,\alpha)$; see Figure \ref{fig:2-3orbitStartingFrom1-2}.  We see that the orbit avoids the peripheral square of $\omegasai{7}{1}$.  By Theorem \ref{thm:AlphaInAFormsOrbitOfSa}, the initial condition $(x^0,\alpha)$ determines an orbit of $\omegasa{7}$. The path traversed by the orbit of $\omegasa{7}$ is exactly the path traversed by the orbit $\mathscr{O}(x^0,\alpha)$.
\end{example}

\section{Concluding remarks}
It is clear from the preceding sections that much work remains to be developed in order to determine a well-defined phase space $(\omegafrac\times S^1)/\sim$ for the yet to be defined fractal billiard flow.  We have discussed several examples of what clearly constitute periodic orbits of $\omegaks$ and $\omegasa{a}$.  Furthermore, for both $\omegaks$ and $\omegat$, we were able to connect two elusive points of each billiard table via suitably chosen nontrivial paths.  These nontrivial paths were determined from suitably chosen sequences of compatible periodic orbits.

\begin{question}
\label{ques:determineAdditionalNontrivialPaths}
Let $F$ be either $\ks$ or $\mathscr{T}$. Suppose that two suitably chosen nontrivial paths converge to two distinct elusive points of $\omegafrac$.  For each of the two elusive points, is it possible to determine another nontrivial path converging to a different elusive point?
\end{question}

If we can answer Question \ref{ques:determineAdditionalNontrivialPaths} in the affirmative (or answer it in the affirmative under specific conditions), will this help us gain insight into how to determine a well-defined phase space for the billiard flow on $\omegafrac$?  An alternate approach, discussed in the concluding remarks of \cite{LapNie3}, entails determining a well-defined fractal translation surface.

Following this line of thought to its logical end, for certain fractal billiard tables (e.g., $\omegat$), is it possible to determine which directions produce recurrent orbits?  More generally, can one prove that, in almost every direction, the billiard flow is ergodic in $\Omega(F)$?

\begin{question}
Regarding a self-similar Sierpinski carpet billiard $\Omega(S_a)$, we have determined a countable set of points from which a periodic billiard orbit can begin.  Can we show that the set of points from which a periodic orbit can begin is in fact uncountable and, furthermore, a set of full \emph{(}Lebesgue\emph{)} measure in the base of the unit square $S_0$?
\end{question}

It is possible to construct a nontrivial line segment of $S_a$ beginning from $(\frac{1}{2},0)$ with slope $\alpha\in \slopesa{a}$, that, when translated to $(0,0)$, no longer lies entirely in $S_{a}$.  However, if we consider the sequence of compatible periodic orbits given by $\compseqixang{0}{(0,0)}{\alpha^0}$, is it possible to determine a well-defined limit? The work of \cite{HuLeTr} may prove useful in further exploring the behavior of a sequence of compatible periodic orbits.  Building on the work of \cite{HuLeTr}, the author of \cite{De} has examined the behavior of nonperiodic orbits in what is an example of what is called a \textit{wind-tree billiard}, and what is also strongly suggestive of a Sierpinski carpet.  Such work may provide insight into examining the behavior of a sequence of compatible dense orbits.

\begin{question}
\label{ques:determiningAWellDefinedBilliardFlow}
In analogy with the prefractal billiard and associated translation surface, can a thorough understanding of the geodesic flow on the limiting \emph{(}and still to be mathematically defined\emph{)} \emph{`}fractal translation surface\emph{'} $\mathcal{S}(F)$ aid us in determining a well-defined billiard flow on $\omegafrac$?
\end{question}

The work in progress in \cite{LapNie4} draws upon the work of Gabriela Weitze-Schmith\"usen \cite{We-Sc} and attempts to answer Question \ref{ques:determiningAWellDefinedBilliardFlow} from an algebraic perspective.

Approaching the problem of determining a well-defined billiard flow on a fractal billiard table from many different points of view may prove useful.  The theories of translation surfaces and rational billiards are intimately tied together and more deeply understood by knowing the structure of what is called the \textit{Veech group} (this being the group studied in, for example, \cite{HuSc,Ve3,Ve4,Vo,We-Sc}). In short, the Veech group of a translation surface $\mathcal{S}(D)$ determined from a rational polygon $D$ is the stabilizer of $\mathcal{S}(D)$.

\begin{question}
Let $\Omega(F)$ be a fractal billiard table, with $F$ being approximated by a suitably chosen sequence of rational polygons $\{F_n\}_{n=0}^\infty$. Is it then possible to construct a Veech group for $\Omega(F)$ \emph{(}or rather, of $\mathcal{S}(F)$\emph{)}, presumably in terms of the Veech groups for the prefractal approximations $\omegafraci{n}$ \emph{(}or rather, of the associated translation surfaces $\mathcal{S}(F_n)$\emph{)}?   Will the knowledge of such a group aid us in determining a well-defined billiard flow on $\Omega(F)$?
\end{question}


\begin{thebibliography}{alpha}

\bibitem[AcST]{AcST} Achdou, Y. Sabot, C., Tchou, N.: Diffusion and propagation problems in some ramified domains with a fractal boundary, \textit{M2AN Math. Model. Numer. Anal.} No. 4, \textbf{40} (2006), 623--652.

\bibitem[Ba]{Ba} Barnsley, M. F.: \textit{SuperFractals}: \textit{Patterns of nature}, Cambridge Univ. Press, New York, 2006.



\bibitem[CheNie]{CheNie} Chen, J. P., Niemeyer, R. G.: Periodic billiard orbits of self-similar Sierpinski carpets, 29 pages, e-print, arXiv:1303.4032v1, 2013. (To appear in the J. of Math. Anal. and Appl.)

\bibitem[De]{De} Delecroix, V.: Divergent directions in some periodic wind-tree models, \textit{J. of Mod. Dyn.} No. 1, \textbf{7} (2013), 1--29.


\bibitem[Du-CaTy]{Du-CaTy} Durand-Cartagena, E., Tyson, J. T.: Rectifiable curves in Sierpi\'{n}ski carpets, \textit{Indiana Univ. Math. J.}  \textbf{60} (2011),  285--310.











\bibitem[Fa]{Fa} Falconer, K. J.: \textit{Fractal Geometry}: \textit{Mathematical foundations and applications}, John Wiley \& Sons, Chichester, 1990. (2nd edition, 2003.)


\bibitem[GaStVo]{GaStVo} Galperin, G., Stepin, A. M., Vorobets, Ya. B.: Periodic billiard trajectories in polygons, \textit{Russian Math. Surveys} No. 3, \textbf{47} (1992), 5--80.



\bibitem[Gut1]{Gut1} Gutkin, E.: Billiards in polygons: Survey of recent results, \textit{J. Stat. Phys.} \textbf{83} (1996), 7--26.

\bibitem[Gut2]{Gut2} Gutkin, E.: Billiards on almost integrable polyhedral surfaces, \textit{Erg. Th. and Dyn. Syst.} \textbf{4} (1984), 569--584.

\bibitem[GutJu1]{GutJu1} Gutkin, E., Judge, C.: The geometry and arithmetic of flat surfaces with applications to polygonal billiards, \textit{Math. Res. Lett.} \textbf{3} (1996), 391--403.

\bibitem[GutJu2]{GutJu2} Gutkin, E., Judge, C.: Affine mappings of flat surfaces: Geometry and arithmetic, \textit{Duke Math. J.} \textbf{103} (2000), 191--213.







\bibitem[HuLeTr]{HuLeTr} Hubert, P., Lelievre, S., Troubetzkoy, S.: The Ehrenfest wind-tree model: periodic directions, recurrence, diffusion, \textit{Journal f\"ur die Reine und Angewandte Mathematik} \textbf{656} (2011), 223--244.

\bibitem[HuSc]{HuSc} Hubert, P., Schmidt, T.: An introduction to Veech surfaces, in: \textit{Handbook of Dynamical Systems}, vol. 1B (A. Katok and B. Hasselblatt, eds.), Elsevier, Amsterdam, 2006, pp. 501--526.

\bibitem[Hut]{Hut} Hutchinson, J. E.: Fractals and self-similarity, \textit{Indiana Univ. Math. J}. \textbf{30} (1981), 713--747.


\bibitem[KaZe]{KaZe} Katok, A., Zemlyakov, A.: Topological transitivity of billiards in polygons, \textit{Math. Notes} \textbf{18} (1975), 760--764.



\bibitem[LapNie1]{LapNie1} Lapidus, M. L., Niemeyer, R. G.: Towards the Koch snowflake fractal billiard---Computer experiments and mathematical conjectures, in: \textit{Gems in Experimental Mathematics} (T. Amdeberhan, L. A. Medina and V. H. Moll, eds.), Contemporary Mathematics, Amer. Math. Soc., Providence, RI, \textbf{517} (2010), pp. 231--263. [E-print: arXiv:math.DS.0912.3948v1, 2009.]

\bibitem[LapNie2]{LapNie2} Lapidus, M. L., Niemeyer, R. G.: Families of periodic orbits of the Koch snowflake fractal billiard, 63 pages, e-print, arXiv:1105.0737v1, 2011.

\bibitem[LapNie3]{LapNie3} Lapidus, M. L., Niemeyer, R. G.: Sequences of compatible periodic hybrid orbits of prefractal Koch snowflake billiards, \textit{Discrete and Continuous Dynamical Systems -- Ser. A}, in press, 2012. [E-print: IHES/M/12/16, 2012; arXiv:1204.3133v1 [math.DS], 2012.]

\bibitem[LapNie4]{LapNie4} Lapidus, M. L., Niemeyer, R. G.: Experimental evidence in support of a fractal law of reflection, in progress, 2013.

\bibitem[LapNie5]{LapNie5} Lapidus, M. L., Niemeyer, R. G.: Veech groups $\Gamma_n$ of the Koch snowflake prefractal translation surfaces $\mathcal{S}(\ks_n)$, in progress, 2012.

\bibitem[LapNie6]{LapNie6} Lapidus, M. L., Niemeyer, R. G.: Sequences of compatible periodic orbits of the $T$-fractal billiard, in progress, 2012.








\bibitem[Mas]{Mas} Masur, H.: Closed trajectories for quadratic differentials with an applications to billiards, \textit{Duke Math. J.} \textbf{53} (1986), 307--314.

\bibitem[MasTa]{MasTa} Masur, H., Tabachnikov, S.: Rational billiards and flat structures, in: \textit{Handbook of Dynamical Systems}, vol. 1A (A. Katok and B. Hasselblatt, eds.), Elsevier, Amsterdam, 2002, pp. 1015--1090.










\bibitem[Sm]{Sm} Smillie, J.: Dynamics of billiard flow in rational polygons, in: \textit{Dynamical Systems}, Encyclopedia of Math. Sciences, vol. 100, Math. Physics 1 (Ya. G. Sinai, ed.), Springer-Verlag, New York, 2000, pp. 360--382.

\bibitem[Ta1]{Ta1} Tabachnikov, S.: \textit{Billiards}, Panoramas et Synth\`eses, Soc. Math. France, Paris, 1995.

\bibitem[Ta2]{Ta2} Tabachnikov, S.: \textit{Geometry and Billiards}, Amer. Math. Soc., Providence, RI, 2005.




\bibitem[Ve1]{Ve1} Veech, W. A.: The billiard in a regular polygon, \textit{Geom. Funct. Anal.} \textbf{2} (1992), 341--379.

\bibitem[Ve2]{Ve2} Veech, W. A.: Flat surfaces, \textit{Amer. J. Math.} \textbf{115} (1993), 589--689.

\bibitem[Ve3]{Ve3} Veech, W.: Teichm\"{u}ller geodesic flow, \textit{Annals of Math}. \textbf{124} (1986), 441--530.

\bibitem[Ve4]{Ve4} Veech, W.: Teichm\"{u}ller curves in modular space, Eisenstein series, and an application to triangular billiards, \textit{Invent. Math.} \textbf{97} (1989), 553--583.

\bibitem[Vo]{Vo} Vorobets, Ya. B.: Plane structures and billiards in rational polygons: The Veech alternative, \textit{Russian Math. Surveys} \textbf{51} (1996), 779--817.

\bibitem[We-Sc]{We-Sc} Weitze-Schmith\"usen, G.: An algorithm for finding the Veech group of an origami, \textit{Experimental Mathematics} No.4, \textbf{13} (2004), 459--472.


\bibitem[Zo]{Zo} Zorich, A.: Flat surfaces, in: \textit{Frontiers in Number Theory, Physics and Geometry} I (P. Cartier, \textit{et al}., eds.), Springer-Verlag, Berlin, 2002, pp. 439--585.
\end{thebibliography}
\end{document}